\documentclass[a4paper, reqno, 11pt]{amsart}
\usepackage[utf8]{inputenc}

\usepackage{amsthm}
\usepackage{amssymb}
\usepackage[american]{babel}
\usepackage{exscale}
\usepackage[mathscr]{euscript}
\usepackage{url}
\usepackage[all,ps,tips,tpic]{xy}
\usepackage{graphicx}
\usepackage{color}

\newenvironment{altenumerate}
   {\begin{list}
      {\textup{(\theenumi)} }
      {\usecounter{enumi}
       \setlength{\labelwidth}{0pt}
       \setlength{\labelsep}{2pt}
       \setlength{\leftmargin}{0pt}
       \setlength{\itemsep}{\the\smallskipamount}
       \renewcommand{\theenumi}{\roman{enumi}}
      }}
   {\end{list}}
\newenvironment{altitemize}
   {\begin{list}
      {$\bullet$ }
      {\setlength{\labelwidth}{0pt}
       \setlength{\labelsep}{2pt}
       \setlength{\leftmargin}{0pt}
       \setlength{\itemsep}{\the\smallskipamount}
      }}
   {\end{list}}

\newtheorem{lem}{Lemma}[section]

\newtheorem{definition}[lem]{Definition}

\newtheorem{thm}[lem]{Theorem}
\newtheorem{prop}[lem]{Proposition}

\theoremstyle{remark}
\newtheorem{rem}[lem]{Remark}

\DeclareMathOperator{\Lie}{Lie}
\DeclareMathOperator{\ord}{ord}

\newcommand{\tr}{\operatorname*{tr}}
\newcommand{\GL}{\mathrm{GL}}

\newdir{ >}{{}*!/-5pt/@{>}}
\SelectTips{cm}{10}

\begin{document}
\title[Deformation spaces of $p$-divisible groups]{The Langlands-Kottwitz method and deformation spaces of $p$-divisible groups}
\author{Peter Scholze}
\begin{abstract}
We extend the results of Kottwitz, \cite{KottwitzPoints}, on points of Shimura varieties over finite fields to cases of bad reduction. The "test function" whose twisted orbital integrals appear in the final expression is defined geometrically using deformation spaces of $p$-divisible groups.
\end{abstract}

\date{\today}
\maketitle
\tableofcontents
\pagebreak

\section{Introduction}

The aim of this article is to extend the results of Kottwitz on the cohomology of Shimura varieties to cases of bad reduction, with the aim of generalizing parts of the results from \cite{ScholzeLLC} to more general Shimura varieties, thereby putting them in their natural context. The results of this paper are used in joint work with S. W. Shin in \cite{ScholzeShin} to prove new results about the cohomology of compact unitary group Shimura varieties at ramified split places, and rederive previously known facts, notably the existence of Galois representations attached to certain automorphic representations of $\GL_n$ over CM fields as in \cite{Shin}.

Let us first briefly recall Kottwitz' results. The basic idea, first due to Langlands in his Antwerp paper, \cite{LanglandsAntwerp}, is to analyze the cohomology of Shimura varieties by computing the alternating sum of the traces of certain Hecke operators twisted by a Frobenius correspondence on the cohomology. In the simplest case, this amounts to counting the number of $\mathbb{F}_{p^r}$-rational points in the special fibre. The key insight is that Honda-Tate theory combined with some group theory allows one to get a purely group-theoretic description of this set. Roughly, it is done in two steps:
\begin{altitemize}
\item Classify the $\mathbb{F}_{p^r}$-isogeny classes.
\item Classify all points within one $\mathbb{F}_{p^r}$-isogeny class.
\end{altitemize}
The result is that the isogeny classes are roughly parametrized by (certain) conjugacy classes in $\mathbf{G}(\mathbb{Q})$, where $\mathbf{G}/\mathbb{Q}$ is the reductive group giving rise to the Shimura variety, and the points within one isogeny class can be described by giving lattices in the \'{e}tale and crystalline cohomology. In the simplest case of the modular curve, this is beautifully explained in unpublished notes of Kottwitz, \cite{KottwitzNotes}, cf. also \cite{Scholze}.

In the work of Kottwitz, \cite{KottwitzPoints}, a similar description of the $\mathbb{F}_{p^r}$-rational points is given for general (compact) PEL Shimura varieties at unramified places. The assumption that one works at an unramified place ensures that the Shimura variety has good reduction and hence that the cohomology of the generic and special fibre agree.

The present work is based on the following observations, which are already present in previous work, cf. e.g. the survey article of Haines, \cite{Haines}.
\begin{altitemize}
\item Even when the PEL data are (mildly, cf. later) ramified, Kottwitz' arguments go through without change to give a description of the $\mathbb{F}_{p^r}$-rational points, for a suitable integral model $\mathcal{M}_{K^p}$, $K^p\subset \mathbf{G}(\mathbb{A}_f^p)$, of the Shimura variety with no level at $p$. This model will however in general not be smooth (not even flat).
\item Even when $\mathcal{M}_{K^p}$ is not flat (but still proper), one can still compute the cohomology of the generic fibre as the cohomology of the special fibre with coefficients in the nearby cycle sheaves, and use the Lefschetz trace formula. In particular, instead of counting the number of fixed points, we have to weight each fixed point with some factor defined in terms of the nearby cycle sheaves.
\item By a theorem of Berkovich, the nearby cycle sheaves depend only on the formal completion of the Shimura variety at the given point, which, in turn, by the theorem of Serre-Tate, depends only on the $p$-divisible group at the given point. This roughly says that the weighting factor corresponding to a point in the isogeny class parametrized by a conjugacy class $\gamma\in \mathbf{G}(\mathbb{Q})$ depends only on the image of $\gamma$ in $\mathbf{G}(\mathbb{Q}_p)$ via some function on $\mathbf{G}(\mathbb{Q}_p)$.\footnote{More precisely, this conjugacy class in $\mathbf{G}(\mathbb{Q}_p)$ will be canonically lifted to a $\sigma$-conjugacy class on $\mathbf{G}(\mathbb{Q}_{p^r})$ defined via the action of the Frobenius operator $F$ on the Dieudonn\'{e} module of the $p$-divisible group, and the function will be defined on $\mathbf{G}(\mathbb{Q}_{p^r})$.} This makes this contribution sufficiently independent of the rest, so that the further formal manipulations allow one to prove a formula close to the one of Kottwitz.
\end{altitemize}

The key new observations of \cite{ScholzeLLC}, with precursors occuring in \cite{Scholze} and \cite{ScholzeGLn}, are the following.
\begin{altitemize}
\item The same method even works if one allows a nontrivial level $K_p$ at $p$, by reinterpreting the cohomology of the Shimura variety of level $K_pK^p$ as the cohomology of the Shimura variety with no level at $p$ (and level $K^p$ away from $p$) with suitable coefficients.
\item For the global applications, it is often enough to know the existence of some $C_c^\infty$ function on $\mathbf{G}(\mathbb{Q}_{p^r})$ that can be used in the trace formula, without any knowledge about its precise values.
\end{altitemize}

The final output is a formula of the form
\begin{equation}\label{e:main-formula}
\mathrm{tr}(\tau\times h f^p|H^{\ast}) = \sum_{\substack{(\gamma_0;\gamma,\delta)\\ \alpha(\gamma_0;\gamma,\delta)=1}} c(\gamma_0;\gamma,\delta) O_{\gamma}(f^p) TO_{\delta\sigma}(\phi_{\tau,h})\ ,
\end{equation}
where $\tau$ is an element of the local Weil group, $h$ is a function in $C_c^\infty(\mathbf{G}(\mathbb{Z}_p))$ (for a suitable integral model of $\mathbf{G}$), and $f^p\in C_c^\infty(\mathbf{G}(\mathbb{A}_f^p))$; we refer to Section \ref{PELShimuraVarieties} for precise assumptions and statements. Let us just mention that the sum basically runs over isogeny classes, which are parametrized in this case not by conjugacy classes in $\mathbf{G}(\mathbb{Q})$, but instead give rise to what we call a Kottwitz triple $(\gamma_0;\gamma,\delta)$ satisfying certain compatibility conditions. The first factor $c(\gamma_0;\gamma,\delta)$ is basically a volume factor as it occurs in trace formulas; the second factor is the orbital integral of $f^p$ and roughly counts lattices in the \'{e}tale cohomology; both of these are well understood. The third factor is the new ingredient; it involves the function $\phi_{\tau,h}\in C_c^\infty(\mathbf{G}(\mathbb{Q}_{p^r}))$ that encodes the weighting factors defined in terms of the nearby cycle sheaves.

In a sense, this puts the whole mystery of the ramification in the bad reduction of the Shimura varieties into certain functions $\phi_{\tau,h}\in C_c^\infty(\mathbf{G}(\mathbb{Q}_{p^r}))$ depending on elements of the local Weil group $\tau$ and a function $h$ on $\mathbf{G}(\mathbb{Z}_p)$.

Let us say a few more words about the construction of this function which is the technical heart of this paper. This also allows us to explain the precise assumptions we put on the PEL data at $p$.

We start by noting that as in the book of Rapoport -- Zink, \cite{RapoportZinkPeriodSpaces}, we consider deformation spaces of $p$-divisible groups, and our main theorem can be read as a theorem that relates the cohomology of Shimura varieties to the cohomology of deformation spaces of $p$-divisible groups. However, these deformation spaces will not be Rapoport-Zink spaces. They occur as formal completions of Rapoport-Zink spaces.

The local PEL data\footnote{We also consider the case of EL data.} that we consider are of the same form as in the book of Rapoport -- Zink, but we will put some additional assumptions. We start with a semisimple $\mathbb{Q}_p$-algebra $B$ with center $F$ such that every simple factor of $B$ is split, i.e. a matrix algebra over a factor of $F$. Let $V$ be a finitely generated left $B$-module. We fix a maximal order $\mathcal{O}_B\subset B$ and a $\mathcal{O}_B$-stable lattice $\Lambda\subset V$. These data give rise to $C=\mathrm{End}_B(V)$ with maximal order $\mathcal{O}_C=\mathrm{End}_{\mathcal{O}_B}(\Lambda)$. Additionally, we have an anti-involution $\ast$ on $B$ which preserves $F$; we let $F_0\subset F$ be the invariants under $\ast$. We make the assumption that $F/F_0$ is unramified. Further, the PEL data consists of a nondegenerate $\ast$-hermitian form $(\, ,\, )$ on $V$, i.e. a nondegenerate alternating form $(\, ,\, )$ on $V$ such that $(bv,w) = (v,b^{\ast}w)$ for all $v, w\in V$, $b\in B$. We require that $\Lambda$ be self-dual with respect to $(\, ,\, )$. This induces an involution $\ast$ on $C$ and $\mathcal{O}_C$. We get the algebraic group $\mathbf{G}/\mathbb{Z}_p$ whose $R$-valued points are given by
\[
\mathbf{G}(R) = \{g\in (R\otimes_{\mathbb{Z}_p} \mathcal{O}_C)^{\times}\mid gg^{\ast}\in R^{\times} \} \ .
\]
The final datum is a conjugacy class $\overline{\mu}$ of cocharacters $\mu: \mathbb{G}_m\longrightarrow \mathbf{G}_{\bar{\mathbb{Q}}_p}$, with field of definition $E$. We assume that after choosing a representative $\mu$ of $\overline{\mu}$, that under the corresponding weight decomposition on $V$, only weights $0$ and $1$ occur, so that $V= V_0\oplus V_1$. Additionally, we assume that the composition
\[
\mathbb{G}_m\buildrel \mu\over\rightarrow \mathbf{G}\rightarrow \mathbb{G}_m
\]
is the identity, the latter morphism denoting the multiplier $g\mapsto gg^{\ast}\in \mathbb{G}_m$.

Moreover, we assume that after extending scalars to $\bar{\mathbb{Q}}_p$, all simple factors of the data $(F,B,\ast,V,(\, ,\, ))$ are of type $A$ or $C$ under the classification of the possible simple factors on page 32 of \cite{RapoportZinkPeriodSpaces}.

Let us mention the restrictions in comparison to \cite{RapoportZinkPeriodSpaces}: We exclude cases of orthogonal type, we assume that the lattice chain is reduced to one selfdual lattice (and translates of it), we assume that $F/F_0$ is unramified, and that $B$ is split.\footnote{But we allow $p=2$, which is excluded in \cite{RapoportZinkPeriodSpaces}.} All of these assumptions are made to avoid additional group-theoretic difficulties; the geometric part of the argument should work without these assumptions, and it would be an interesting problem to extend this method beyond the cases considered here.

Let $\mathcal{D}=(B,F,V,\ldots)$ denote the PEL data. One gets a natural notion of $p$-divisible group with $\mathcal{D}$-structure over any scheme $S$ on which $p$ is locally nilpotent, cf. Definition \ref{DefPDivGroupDStr}. Over $\mathbb{F}_{p^r}$, one can describe them via Dieudonn\'{e} theory by elements $\delta\in \mathbf{G}(\mathbb{Q}_{p^r})$ up to $\sigma$-conjugation by $\mathbf{G}(\mathbb{Z}_{p^r})$, cf. Proposition \ref{ParamPDivGroupDStrPerfField}. Let $\overline{\underline{H}} = \overline{\underline{H}}_{\delta}$ be the $p$-divisible group with $\mathcal{D}$-structure associated to some $\delta$. Then we look at its universal deformation which lives over a complete noetherian local ring $R_{\overline{\underline{H}}}$, which in turn gives rise to a rigid-analytic variety $X_{\overline{\underline{H}}}$. Moreover, for any $K\subset \mathbf{G}(\mathbb{Z}_p)$, we get a natural finite \'{e}tale cover $X_{\overline{\underline{H}},K}/X_{\overline{\underline{H}}}$. Let $X_{\overline{\underline{H}},K,\overline{\eta}}$ be the base-change to $\mathbb{C}_p$.

The basic idea is to define
\[
\phi_{\tau,h}(\delta) = \tr(\tau\times h | H^\ast(X_{\overline{\underline{H}},K,\overline{\eta}},\mathbb{Q}_\ell))\ ,
\]
for any $K\subset \mathbf{G}(\mathbb{Z}_p)$ such that $h$ is biinvariant under $K$.\footnote{Here, we use Huber's definition of \'{e}tale cohomology groups for rigid-analytic varieties.} Note that the nearby cycles can be identified with the cohomology of the generic fibre, thus this expression really is the local term of the Lefschetz trace formula that occurs as the weighting factor.

At this point, we can explain the main technical obstacle: It is not in general known that the cohomology groups of $X_{\overline{\underline{H}},K}$ satisfy suitable finiteness statements. This would follow if one could prove that these spaces are algebraizable in a suitable sense, cf. Theorem \ref{AlgebraizationControlled}. In the special case considered in our previous work \cite{ScholzeLLC}, we could deduce this from Faltings's theory of group schemes with strict $\mathcal{O}$-action, \cite{FaltingsStrictOAction}, and Artin's algebraization theorem. The same proof works for unramified PEL data (and $p\neq 2$), by using a theorem of Wedhorn, \cite{WedhornOortStrata}, Theorem 2.8, which describes the deformation space as a versal deformation space of a truncated $p$-divisible group, cf. Proposition \ref{AlgUnram}. However, it seems difficult to us to make this strategy of proving algebraicity work outside the case where the moduli problem without level structure is smooth. In joint work with S. W. Shin, \cite{ScholzeShin}, we will use a global argument to prove that in all EL cases, the deformation spaces are algebraizable by showing that they occur in a suitable Shimura variety. However, this argument will already make use of results proven here.

The idea to get around this difficulty is the observation that in the end, we are only interested in $\phi_{\tau,h}(\delta)$ if the $p$-divisible group associated with $\delta$ comes from some point in the Shimura variety. In that case, there tautologically is an algebraization, and hence the finiteness results that one needs hold true.

We employ this idea by defining a notion of rigid-analytic varieties with "controlled cohomology"; this implies that the cohomology is the same as the cohomology of a quasicompact admissible open subset. For quasicompact rigid-analytic varieties (satisfying slight technical extra assumptions), finiteness results are known, and hence the cohomology of any rigid-analytic variety with "controlled cohomology" satisfies suitable finiteness statements. It also implies all other results that one would like to know, like independence of $\ell$, etc. .

This leads to a well-defined function $\phi_{\tau,h}$ on $\mathbf{G}(\mathbb{Q}_{p^r})$, and the final statement one needs is that this is a locally constant function. This is based on the result that the automorphism group of $\underline{\overline{H}}$ acts smoothly on the cohomology (which follows from results of Berkovich, but which we deduce here from results of Huber).

Let us add some remarks about previous works using the Langlands-Kottwitz method in cases of bad reduction. The first case considered was that of parahoric level structures. In that case, there are natural integral models for the Shimura variety with parahoric level structure (which are not used here), and one can adapt Kottwitz's arguments to describe the points in the special fibre of these integral models of parahoric level. A similar function $\phi^\prime$ can be constructed, and a conjecture of Kottwitz states that this function lies in the center of the Iwahori-Hecke algebra, and can be explicitly described via the Bernstein isomorphism. This led Beilinson to a general conjecture in the geometric Langlands program that one can construct central elements in Iwahori-Hecke algebras via a nearby-cycles construction, which was proved by Gaitsgory, \cite{GaitsgoryCentralNearby}. His argument was adapted by Haines and Ng\^o to prove Kottwitz' original conjecture, \cite{HainesNgo}. We note that this conjecture of Kottwitz leads to a precise formula for the spectral contribution of the test function $\phi^\prime$ that becomes necessary in applications of formula \eqref{e:main-formula}, e.g. the computation of the semisimple local Hasse-Weil zeta function as in \cite{Haines}, Theorem 11.7. Outside the parahoric case, we mention an article of Haines and Rapoport, \cite{HainesRapoport}, which considers the unipotent radical of an Iwahori subgroup for a specific Shimura variety. Again, they find that there is a canonical test function that lies in the center of a Hecke algebra, and they identify it explicitly, which leads to a description of the Hasse-Weil zeta function.

This led to the expectation that for general level structures, it should be possible to define canonical test functions in the center of an appropriate Hecke algebra. In the paper \cite{Scholze} dealing with the case of modular curve, it was verified that functions in the center of a Hecke algebra exist which can be used to make \eqref{e:main-formula} true. These functions were made explicit, \cite{Scholze}, Section 14, but they were constructed artificially, and no relation to the geometry could be found. It then became clear, \cite{ScholzeGLn}, \cite{ScholzeLLC}, that one can define canonical test functions coming from the geometry, but that these do not lie in the center of a Hecke algebra. In \cite{ScholzeGLn}, it was verified formally that also functions in the center of a Hecke algebra exist, but no geometric interpretation could be given.

Nonetheless, there is an analogue of Kottwitz' conjecture, as formulated in joint work with S. W. Shin, \cite{ScholzeShin}. However, it gives less information: Our conjecture only determines the twisted orbital integrals of $\phi_{\tau,h}$, whereas Kottwitz' conjecture determined the function itself. Even in the case that we take $h$ as the idempotent associated to a parahoric level structure, our function $\phi_{\tau,h}$ behaves worse than the function $\phi^\prime$ in that it does not lie in the center of the Iwahori-Hecke algebra. This is due to the fact that different integral models of the Shimura variety are used to define them. One expects that their twisted integral orbitals agree, which is however not clear from the definition.

Finally, we give an overview of the content. In Section 2, we collect some general results about the cohomology of rigid-analytic varieties as proved by Huber. In Section 3, we define the deformation spaces of $p$-divisible groups, which we use in Section 4 to define the test functions $\phi_{\tau,h}$. In Section 5, we state our main theorem calculating the trace of certain Hecke operators twisted with an element of the local Weil group on the cohomology of Shimura varieties. This theorem is proved in Sections 6 and 7; in Section 6, we give a description of the fixed points of Hecke correspondences, and in Section 7, we apply the Lefschetz trace formula.

{\bf Acknowledgments.} It is a pleasure to thank my advisor M. Rapoport for everything he taught me, and his continuous help with problems of any sort. Moreover, I thank T. Haines and R. Kottwitz for helpful discussions; this paper obviously owes a lot to their work. This paper was written while the author was a Clay Research Fellow.

\section{Etale cohomology of rigid-analytic varieties}

In this section, we recall some facts about the \'{e}tale cohomology of rigid-analytic varieties that will be used.

In the following, we fix a complete discrete valuation field $k$ of characteristic $0$, with valuation subring $k^\circ$ and maximal ideal $k^{\circ\circ}$. Moreover, we fix a prime $\ell$ prime to the characteristic of the residue class field $\kappa = k^\circ/k^{\circ\circ}$.

We have the following properties, collected from the work of Huber, \cite{Huber}.

\begin{thm} Let $X$ and $Y$ be separated smooth rigid-analytic varieties over $k$ of dimension $n$, resp. $m$.
\begin{altenumerate}
\item[{\rm (i)}] If $X$ is quasicompact, then the cohomology groups $H^i(X\otimes_k \hat{\bar{k}},\mathbb{Z}/\ell^j\mathbb{Z})$ and $H^i_c(X\otimes_k \hat{\bar{k}}, \mathbb{Z}/\ell^j\mathbb{Z})$ are finite for all $i$ and $j$.
\item[{\rm (ii)}] Assume that $X$ is taut (cf. Definition 0.4.7 of \cite{Huber}). Then the cohomology groups $H^i(X\otimes_k \hat{\bar{k}},\mathbb{Z}/\ell^j\mathbb{Z})$ and $H^i_c(X\otimes_k \hat{\bar{k}}, \mathbb{Z}/\ell^j\mathbb{Z})$ are zero for $i>2n$. Moreover, we have an isomorphism
\[
H^i(X\otimes_k \hat{\bar{k}},\mathbb{Z}/\ell^j\mathbb{Z})\cong H^{2n-i}_c(X\otimes_k \hat{\bar{k}},\mathbb{Z}/\ell^j\mathbb{Z})^{\vee}(-n)
\]
for $i=0,\ldots,2n$, where $(-n)$ denotes a Tate twist.
\item[{\rm (iii)}] Assume that $X$ and $Y$ are taut. Then there is a K\"unneth-formula isomorphism
\[
R\Gamma_c(X\otimes_k \hat{\bar{k}},\mathbb{Z}/\ell^j \mathbb{Z})\buildrel{\mathbb{L}}\over \otimes_{\mathbb{Z}/\ell^j \mathbb{Z}} R\Gamma_c(Y\otimes_k \hat{\bar{k}},\mathbb{Z}/\ell^j \mathbb{Z})\cong R\Gamma_c((X\times Y)\otimes_k \hat{\bar{k}},\mathbb{Z}/\ell^j\mathbb{Z})\ .
\]
\end{altenumerate}
\end{thm}

\begin{proof} Part (i) follows from Proposition 0.5.3 and Proposition 0.5.4 of \cite{Huber}. Part (ii) follows from Corollary 0.5.8 and Corollary 0.6.3 of \cite{Huber}. Part (iii) follows as usual from the properties 0.4.5 a) - d) of \cite{Huber}, noting that they hold true if $f$ is separated and taut, cf. \cite{Huber}, p. 19.
\end{proof}

Now assume that $X$ is a separated taut smooth rigid-analytic variety over $k$ of dimension $n$, and assume that we are given a sequence of quasicompact admissible open subspaces $X_0\subset X_1\subset \ldots \subset X$, such that $X=\bigcup_j X_j$.

\begin{definition} We say that $X$ has \emph{controlled cohomology} (with respect to $X_0,X_1,\ldots$) if for all large $j$, the map
\[
H^i_c(X_j\otimes_k \hat{\bar{k}},\mathbb{Z}/\ell\mathbb{Z})\rightarrow H^i_c(X\otimes_k \hat{\bar{k}},\mathbb{Z}/\ell\mathbb{Z})
\]
is an isomorphism for all $i$.
\end{definition}

\begin{prop} Assume that $X$ is a separated taut smooth rigid-analytic variety over $k$ of dimension $n$ having controlled cohomology with respect to a sequence of subsets $X_0,X_1,\ldots$ as above. Then the groups $H^i(X\otimes_k \hat{\bar{k}},\mathbb{Z}/\ell^j\mathbb{Z})$ and $H^i_c(X\otimes_k \hat{\bar{k}}, \mathbb{Z}/\ell^j\mathbb{Z})$ are finite for all $i$ and $j$, and there is a perfect pairing
\[
H^i_{\mathrm{et}}(X\otimes_k \hat{\bar{k}},\mathbb{Z}/\ell^j\mathbb{Z})\otimes H^{2n-i}_{c,\mathrm{et}}(X\otimes_k \hat{\bar{k}},\mathbb{Z}/\ell^j\mathbb{Z})\rightarrow \mathbb{Z}/\ell^j \mathbb{Z}(-n)
\]
for all $i\in \mathbb{Z}$; in particular the cohomology groups vanish for $i>2n$.
\end{prop}

\begin{proof} Immediate.
\end{proof}

\begin{prop}\label{ProductControlled} Let $X$ and $Y$ be seperated taut smooth rigid-analytic varieties over $k$ of dimension $n$, resp. $m$, and consider quasicompact admissible open subsets $X_0\subset X_1\subset \ldots \subset X$ and $Y_0\subset Y_1\subset \ldots\subset Y$ exhausting $X$, resp. $Y$. Then $X\times Y$ with respect to the sequence of quasicompact admissible open subsets $X_0\times Y_0\subset X_1\times Y_1\subset \ldots\subset X\times Y$ has controlled cohomology if and only if both $X$ and $Y$ have controlled cohomology (with respect to the given $X_j$, resp. $Y_j$).
\end{prop}

\begin{proof} Use the K\"unneth formula.
\end{proof}

\begin{prop}\label{GaloisCoverCohom} Let $X$ be a separated taut smooth rigid-analytic variety over $k$ of dimension $n$ with $X_0\subset X_1\subset \ldots \subset X$ as before, and let $f: Y\rightarrow X$ be a finite \'{e}tale Galois cover with Galois group $G$. Let $Y_j = f^{-1}(X_j)$; then $Y_0\subset Y_1\subset \ldots \subset Y$ are quasicompact admissible open subsets exhausting $Y$, and $Y_j\rightarrow X_j$ is a finite \'{e}tale Galois cover with Galois group $G$ for all $j$. Assume that the order of $G$ is prime to $\ell$. Then if $Y$ has controlled cohomology (w.r.t. the $Y_j$), then $X$ has controlled cohomology (w.r.t. the $X_j$), and
\[
H_c^i(X\otimes_k \hat{\bar{k}},\mathbb{Z}/\ell \mathbb{Z}) = H_c^i(Y\otimes_k \hat{\bar{k}}, \mathbb{Z}/\ell\mathbb{Z})^G
\]
for all $i\in \mathbb{Z}$.
\end{prop}

\begin{proof} First, note that $Rf_! \mathbb{Z}/\ell\mathbb{Z}=f_! \mathbb{Z}/\ell\mathbb{Z}=f_\ast \mathbb{Z}/\ell \mathbb{Z}=Rf_\ast \mathbb{Z}/\ell \mathbb{Z}$ is locally constant on $X_{\mathrm{et}}$ by Corollary 0.5.6 of \cite{Huber}. We may decompose it into a direct sum according to the irreducible representations of $G$ over $\mathbb{F}_{\ell}$, using that $\ell$ is prime to the order of $G$. The $G$-invariant part is just $\mathbb{Z}/\ell \mathbb{Z}$ (via the adjunction morphism $\mathbb{Z}/\ell \mathbb{Z}\rightarrow f_\ast \mathbb{Z}/\ell \mathbb{Z}$). This shows that for all $j$, we have
\[
H_c^i(X_j\otimes_k \hat{\bar{k}},\mathbb{Z}/\ell \mathbb{Z}) = H_c^i(Y_j\otimes_k \hat{\bar{k}}, \mathbb{Z}/\ell\mathbb{Z})^G
\]
as well as
\[
H_c^i(X\otimes_k \hat{\bar{k}},\mathbb{Z}/\ell \mathbb{Z}) = H_c^i(Y\otimes_k \hat{\bar{k}}, \mathbb{Z}/\ell\mathbb{Z})^G\ .
\]
The conclusion follows.
\end{proof}

We want to extend this discussion to $\ell$-adic coefficients. In \cite{HuberEllAdic}, Huber defines compactly supported cohomology with $\ell$-adic coefficients. We set
\[
H_c^\ast(X\otimes_k \hat{\bar{k}},\mathbb{Q}_\ell) = H_c^\ast(X\otimes_k \hat{\bar{k}},\mathbb{Z}_\ell)\otimes_{\mathbb{Z}_\ell} \mathbb{Q}_\ell
\]
and
\[
H^i(X\otimes_k \hat{\bar{k}},\mathbb{Q}_\ell) = \mathrm{Hom}(H_c^{2n-i}(X\otimes_k \hat{\bar{k}},\mathbb{Q}_\ell),\mathbb{Q}_\ell)
\]
in case $X$ is smooth of dimension $n$.

Now assume that $X$ is separated, taut and smooth, and has controlled cohomology with respect to $X_j\subset X$. Then Proposition 2.1 iv) and Theorem 3.1 of \cite{HuberEllAdic} imply that
\[
H_c^\ast(X\otimes_k \hat{\bar{k}},\mathbb{Z}_\ell)
\]
is the inverse limit of $H_c^\ast(X\otimes_k \hat{\bar{k}},\mathbb{Z}/\ell^m\mathbb{Z})$, and this inverse system is AR-$\ell$-adic. In particular $H_c^\ast(X\otimes_k \hat{\bar{k}},\mathbb{Q}_\ell)$ and $H^\ast(X\otimes_k \hat{\bar{k}})$ are finite-dimensional $\mathbb{Q}_\ell$-vector spaces vanishing outside the range $0\leq i\leq 2\dim X$, and satisfy a K\"unneth formula, also for ordinary cohomology:
\[
H^i((X\times Y)\otimes_k \hat{\bar{k}},\mathbb{Q}_\ell ) = \bigoplus_j H^j(X\otimes_k \hat{\bar{k}},\mathbb{Q}_\ell)\otimes_{\mathbb{Q}_\ell} H^{i-j}(Y\otimes_k \hat{\bar{k}},\mathbb{Q}_\ell)\ .
\]
Moreover, in the situation of Proposition \ref{GaloisCoverCohom}, we have
\[
H^i(X\otimes_k \hat{\bar{k}},\mathbb{Q}_\ell ) = H^i(Y\otimes_k \hat{\bar{k}},\mathbb{Q}_\ell )^G\ .
\]

Finally, we have the following result that will imply smoothness of certain group actions on the cohomology of the spaces considered.

\begin{prop}\label{Continuity} Let $X$ be a quasicompact separated taut smooth rigid-analytic variety over $k$. Then there is an admissible open neighborhood $U\subset X\times X$ of the diagonal $X\subset X\times X$ such that for all $f: X\rightarrow X$ for which the graph $\Gamma_f\subset X\times X$ is contained in $U$, the induced morphism
\[
f^\ast: H^i(X\otimes_k \hat{\bar{k}},\mathbb{Q}_\ell)\rightarrow H^i(X\otimes_k \hat{\bar{k}},\mathbb{Q}_\ell)
\]
is the identity for all $i\in \mathbb{Z}$.
\end{prop}

\begin{proof} It is enough to ensure that $H_c^i(X\otimes_k \hat{\bar{k}},\mathbb{Z}/\ell\mathbb{Z})$ maps isomorphically to $H_c^i(U\otimes_k \hat{\bar{k}},\mathbb{Z}/\ell\mathbb{Z})$ for all $i$, for then the same is true for $\mathbb{Q}_\ell$-cohomology, and we have a commutative diagram
\[\xymatrix{
& H^i(X\otimes_k \hat{\bar{k}},\mathbb{Q}_\ell) \ar[dr]^{\cong} & \\
H^i(X\otimes_k \hat{\bar{k}},\mathbb{Q}_\ell) \ar[ur]^{\cong} \ar[dr]^{\cong} & H^i(U\otimes_k \hat{\bar{k}},\mathbb{Q}_\ell) \ar[u]^{\cong} \ar[d] & H^i(X\otimes_k \hat{\bar{k}},\mathbb{Q}_\ell) \\
& H^i(\Gamma_f\otimes_k \hat{\bar{k}},\mathbb{Q}_\ell) \ar[ur]^{\cong} &
}\]
The existence of $U$ with this property follows from Theorem 2.9 of \cite{HuberFinitenessCompact}.
\end{proof}

Now we put ourself into the following situation. Let $R$ be a complete noetherian local $k^\circ$-algebra with residue field $\kappa$ and maximal ideal $\mathfrak{m}$. Then to $R$, one can associate a rigid-analytic variety $X$ over $k$ as in \cite{RapoportZinkPeriodSpaces}, 5.5, cf. also \cite{BerthelotCohoRig}, 0.2.6.

Let us recall the construction, in a slightly more general situation. Let us only assume that $R$ is a complete noetherian semilocal $k^\circ$-algebra whose residue fields are finite extensions of $\kappa$, and let $\mathfrak{m}\subset R$ be the set of topologically nilpotent elements. Let $f_1,\ldots,f_m$ be generators of $\mathfrak{m}$, and let $\varpi\in k^{\circ\circ}$ be a uniformizer. Then for any $j\geq 1$, consider the algebra
\[
R_j = R\langle T_1,\ldots,T_m\rangle / (f_1^j - \varpi T_1,\ldots,f_m^j - \varpi T_m)\ .
\]
One checks that $R_j\otimes k$ is a Tate algebra over $k$ in the sense of rigid geometry, and we can define $X_j = \mathrm{Sp}\ R_j$, a rigid-analytic variety over $k$. Moreover, there are obvious transition maps $R_j\rightarrow R_{j^\prime}$ for $j\geq j^\prime$, inducing maps $X_{j^\prime}\rightarrow X_j$ for $j^\prime\leq j$; these are admissible open embeddings. Finally, one defines $X=\bigcup_j X_j$, with the $X_j$ as admissible open subsets.

\begin{lem} For any $j\geq 1$, the subset $X_j\subset X$ is the subset of all $x\in X$ where for all $f\in \mathfrak{m}$, the inequality $|f(x)|\leq |\varpi|^{1/j}$ holds true. The rigid-analytic variety $X$ is separated and partially proper (cf. Definition 0.4.2 in \cite{Huber}), in particular taut.
\end{lem}

\begin{proof} This follows easily from the construction of $X$.
\end{proof}

Now assume moreover that $X$ is smooth of pure dimension $n$. Let $f: Y\rightarrow X$ be a finite \'{e}tale morphism. Taking $Y_j = f^{-1}(X_j)$, we similarly exhaust $Y$ by the quasicompact admissible open subsets $Y_j\subset Y$. In the following, we simply say that $X$ or $Y$ has controlled cohomology, this family of quasicompact open subsets being understood.

The results of \cite{HuberFiniteness} and \cite{HuberFinitenessCompact} imply the following theorem.

\begin{thm}\label{AlgebraizationControlled} Assume that there is a separated scheme $\mathfrak{X}$ of finite type over $k^\circ$ with smooth generic fibre and with a finite \'{e}tale cover $g: \mathfrak{Y}\rightarrow \mathfrak{X}\otimes_{k^\circ} k$. Let $x\in \mathfrak{X}(\kappa)$ such that $R$ is the completed local ring of $\mathfrak{X}$ at $x$, identifying the tubular neighborhood of $x$ in the generic fibre $\mathfrak{X}^{\mathrm{rig}}$ of $\mathfrak{X}$ with $X$, and such that there is a fibre product diagram
\[\xymatrix{
Y\ar[r]\ar@{^(->}[d] & X\ar@{^(->}[d]\\
\mathfrak{Y}^{\mathrm{rig}}\ar[r] & \mathfrak{X}^{\mathrm{rig}}
}\]
Then $X$ and $Y$ have controlled cohomology. Moreover, fixing a geometric point $\overline{x}$ above $x$, there is a canonical $\mathrm{Gal}(\bar{k}/k)$-equivariant isomorphism
\[
(R^i\psi g_\ast \mathbb{Q}_\ell)_{\overline{x}}\cong H^i(Y\otimes_k \hat{\bar{k}},\mathbb{Q}_\ell)\ .
\]
In particular, if $\mathfrak{Y}/\mathfrak{X}$ is Galois with Galois group $G$, so is $Y/X$ and
\[
H^i(X\otimes_k \hat{\bar{k}},\mathbb{Q}_\ell)\cong H^i(Y\otimes_k \hat{\bar{k}},\mathbb{Q}_\ell)^G\ .
\]
\end{thm}

\begin{proof} By Theorem 2.9 of \cite{HuberFinitenessCompact}, both $X$ and $Y$ have controlled cohomology. Proposition 3.15 of \cite{HuberFiniteness} gives the description of the nearby cycles for torsion coefficients. In order to pass to $\ell$-adic coefficients, we use Proposition 5.9.4 from \cite{Fargues}, which shows that Proposition 5.9.2 of \cite{Fargues} applies, giving that the definition of $H^\ast(X\otimes_k \hat{\bar{k}},\mathbb{Q}_\ell)$ agrees with the naive definition as an inverse limit of $H^\ast(X\otimes_k \hat{\bar{k}},\mathbb{Z}/\ell^m\mathbb{Z})$ tensored with $\mathbb{Q}_\ell$ (and the same statement for $Y$ or with coefficients $g_\ast \mathbb{Z}/\ell^m\mathbb{Z}$). The last statement follows from $(g_\ast \mathbb{Q}_\ell)^G = \mathbb{Q}_\ell$ in the algebraic context.
\end{proof}

Finally, we have some lemmas that guarantee that one can apply Proposition \ref{Continuity}.

\begin{lem} Let $R$ be a topological $k^\circ$-algebra, topologically of finite type, and let $X$ be its generic fibre as a rigid-analytic variety over $k$. Then for any open neighborhood $U\subset X\times X$ containing the diagonal $X\subset X\times X$, there is some $a\in k^{\circ\circ}$ such that for all automorphisms $f$ of $R$ that act trivially on $R/a$, the graph $\Gamma_f\subset X\times X$ of $f$ acting on $X$ is contained in $U$.
\end{lem}

\begin{proof} Let $\pi_i: X\times X\rightarrow X$, $i=1,2$, be the two projections. Then a cofinal system of admissible open neighborhoods of the diagonal is given by finite intersections of admissible open neighborhoods of the form $|\pi_1^\ast(r) - \pi_2^\ast(r)|\leq |\varpi|^m$ for $r\in R$, $m\geq 0$. Then any automorphism $f$ that is trivial modulo $\varpi^m$ for all $m$ occuring in these expressions has the desired property.
\end{proof}

\begin{lem}\label{ContinuityHelp} Let $R$ be a complete noetherian semilocal $k^\circ$-algebra whose residue fields are finite extensions of $\kappa$. Let $X$ be the corresponding generic fibre as a rigid-analytic variety over $k$. Recall that $X=\bigcup_j X_j$ is naturally the union of quasicompact admissible open subsets $X_j\subset X$. Then for any $j\geq 1$ and any admissible open neighborhood $U\subset X_j\times X_j$ of the diagonal $X_j\subset X_j\times X_j$, there exists an open ideal $I\subset R$ such that for all automorphisms $f$ of $R$ that act trivially on $R/I$, the graph $\Gamma_f\subset X_j\times X_j$ of the action of $f$ on $X_j$ is contained in $U$.
\end{lem}

\begin{proof} The subset $X_j\subset X$ is affinoid, equal to the generic fibre of $\mathrm{Spf}\ R_j$ for some $R_j$ as in the previous lemma. Now apply the previous lemma.
\end{proof}

\section{Deformation spaces of $p$-divisible groups}\label{DefSpaces}

We will consider certain deformation spaces of $p$-divisible groups with extra structure. This extra structure will be given by endomorphism and level structures (EL case), or polarization, endomorphism and level structures (PEL case). Our conventions are the following, cf. \cite{RapoportZinkPeriodSpaces}, 1.38 and Definition 3.18:
\\

{\bf EL case}. Let $B$ be a semisimple $\mathbb{Q}_p$-algebra with center $F$ such that every simple factor of $B$ is split, i.e. a matrix algebra over a factor of $F$. Let $V$ be a finitely generated left $B$-module. We fix a maximal order $\mathcal{O}_B\subset B$ and a $\mathcal{O}_B$-stable lattice $\Lambda\subset V$. These data give rise to $C=\mathrm{End}_B(V)$ with maximal order $\mathcal{O}_C=\mathrm{End}_{\mathcal{O}_B}(\Lambda)$ and the algebraic group $\mathbf{G}/\mathbb{Z}_p$ whose group of $R$-valued points is given by
\[
\mathbf{G}(R) = (R\otimes_{\mathbb{Z}_p} \mathcal{O}_C)^{\times}
\]
for any $\mathbb{Z}_p$-algebra $R$. The final datum is a conjugacy class $\overline{\mu}$ of cocharacters $\mu: \mathbb{G}_m\longrightarrow \mathbf{G}_{\bar{\mathbb{Q}}_p}$. The field of definition of $\overline{\mu}$ is called $E$, the reflex field. Fix a representative $\mu$ of $\overline{\mu}$ over $\bar{\mathbb{Q}}_p$; this gives rise to a decomposition of $V$ into weight spaces. We make the assumption that only weights $0$ and $1$ occur in the decomposition, i.e. $V=V_0\oplus V_1$. The isomorphism class of the $B$-module $V_0$ (and that of $V_1$) is defined over $E$.
\\

{\bf PEL case}. Again, we have $B$, $\mathcal{O}_B$, $F$, $V$, $\Lambda$, $C$ and $\mathcal{O}_C$ as in the PEL case. Additionally, we have an anti-involution $\ast$ on $B$ which preserves $F$; we let $F_0\subset F$ be the invariants under $\ast$. We make the assumption that $F/F_0$ is unramified. Further, the PEL data consists of a nondegenerate $\ast$-hermitian form $(\, ,\, )$ on $V$, i.e. a nondegenerate alternating form $(\, ,\, )$ on $V$ such that $(bv,w) = (v,b^{\ast}w)$ for all $v, w\in V$, $b\in B$. We require that $\Lambda$ be self-dual with respect to $(\, ,\, )$. This induces an involution $\ast$ on $C$ and $\mathcal{O}_C$. We get the algebraic group $\mathbf{G}/\mathbb{Z}_p$ given by
\[
\mathbf{G}(R) = \{g\in (R\otimes_{\mathbb{Z}_p} \mathcal{O}_C)^{\times}\mid gg^{\ast}\in R^{\times} \} \ .
\]
The final datum is again a conjugacy class $\overline{\mu}$ of cocharacters $\mu: \mathbb{G}_m\longrightarrow \mathbf{G}_{\bar{\mathbb{Q}}_p}$, with field of definition $E$. Again, we assume that after choosing a representative $\mu$ of $\overline{\mu}$, that under the corresponding weight decomposition on $V$, only weights $0$ and $1$ occur, so that $V= V_0\oplus V_1$, where the isomorphism class of the $B$-module $V_0$ is defined over $E$ again. Additionally, we assume that the composition
\[
\mathbb{G}_m\buildrel \mu\over\rightarrow \mathbf{G}\rightarrow \mathbb{G}_m
\]
is the identity, the latter morphism denoting the multiplier $g\mapsto gg^{\ast}\in \mathbb{G}_m$.

Moreover, we assume that after extending scalars to $\bar{\mathbb{Q}}_p$, all simple factors of the data $(F,B,\ast,V,(\, ,\, ))$ are of type $A$ or $C$ under the classification of the possible simple factors on page 32 of \cite{RapoportZinkPeriodSpaces}.
\\

There are certain cases of PEL type which are closely related to cases of EL type. They are needed when one wants to embed the EL case into a global situation, as Shimura varieties are associated to PEL data.
\\

{\bf Quasi-EL case}. Assume given data $\mathcal{D}$ of PEL type such that $F$ decomposes as $F=F_0\times F_0$, with $\ast$ acting by $(x,y)^\ast = (y,x)$. We let
\[
\mathcal{D}_0 = (B_0 = B\otimes_F F_0, \mathcal{O}_{B_0} = \mathcal{O}_B\otimes_{\mathcal{O}_F} \mathcal{O}_{F_0}, V_0 = V\otimes_F F_0, \Lambda_0 = \Lambda\otimes_{\mathcal{O}_F} \mathcal{O}_{F_0},\overline{\mu}_0)\ ,
\]
where $\overline{\mu}_0$ is defined as the first factor of $\overline{\mu}$ under the isomorphism
\[
\mathbf{G} = \mathbf{G}_0\times \mathbb{G}_m
\]
sending an element $g$ acting on $\Lambda = \Lambda_0 \oplus \Lambda_0^\ast$ to its restriction to $\Lambda_0$, and the multiplier $gg^\ast\in \mathbb{G}_m$. This defines data $\mathcal{D}_0$ of EL type. Note that the second component of $\overline{\mu}$ is a morphism $\mathbb{G}_m\rightarrow \mathbb{G}_m$ which by assumption is the identity. In this way, there is a bijection between data of EL type and data of quasi-EL type (with fixed decomposition $F=F_0\times F_0$).
\\

In particular, we note that our assumptions imply that the reductive group $\mathbf{G}_{\mathbb{Q}_p}$ is connected and quasisplit, with simply connected derived group. In comparison with the assumptions of \cite{RapoportZinkPeriodSpaces}, we make the additional assumptions that $B$ is split over $F$, that $F/F_0$ is unramified, that the lattice chain is reduced to a single selfdual lattice $\Lambda\subset V$ (and its translates), and we exclude the cases of orthogonal type. It would be interesting to investigate whether our method can be adapted to more general cases.

The following lemma is crucial at a large number of places, and justifies most of our assumptions.

\begin{lem}\label{CrucialLemma} Let data of PEL type be given. Let $(V^{\prime}, (\, ,\, )^{\prime})$ be a $\ast$-hermitian $B$-module such that $V^{\prime}\cong V$ as $B$-modules. Assume that there is a self-dual $\mathcal{O}_B$-stable lattice $\Lambda^{\prime}\subset V^{\prime}$. Then there is an isomorphism of $\ast$-hermitian $B$-modules $V\cong V^{\prime}$ carrying $\Lambda$ into $\Lambda^{\prime}$.
\end{lem}

\begin{rem} We require the isomorphism to preserve the form, not just up to a scalar. Hence this is slightly stronger than Lemma 7.2 in \cite{KottwitzPoints}. We note that the proposition stays true after tensoring everything with $\mathbb{Q}_{p^r}$ over $\mathbb{Q}_p$.
\end{rem}

\begin{proof} We may assume that $F_0$ is a field. Let $\varpi\in F_0$ be a uniformizer of $F_0$. The group of $\mathcal{O}_B$-linear hermitian isomorphisms of $\Lambda$ is an unramified group over the ring of integers of $F_0$. It is also connected by our assumption that all simple factors of the PEL data are of type $A$ or $C$. By Lang's lemma, there is an isomorphism of $\ast$-hermitian $\mathcal{O}_B$-modules $\Lambda/\varpi\Lambda\cong \Lambda^{\prime} / \varpi\Lambda^{\prime}$. As in \cite{KottwitzPoints}, proof of Lemma 7.2, one checks that this isomorphism lifts.
\end{proof}

Now let $\mathcal{D}$ be data of type EL or PEL. We will consider the following type of $p$-divisible groups with extra structure, cf. Definition 3.21 of \cite{RapoportZinkPeriodSpaces}.

\begin{definition}\label{DefPDivGroupDStr} Let $S$ be a scheme over $\mathcal{O}_E$ on which $p$ is locally nilpotent. A $p$-divisible group with $\mathcal{D}$-structure over $S$ is given by a pair $\underline{H} = (H,\iota)$ (resp. a quadruple $\underline{H}=(H,\iota,\lambda,\mathbb{L})$ in the PEL case), consisting of
\begin{altitemize}
\item a $p$-divisible group $H$ over $S$
\item a homomorphism $\iota:\mathcal{O}_B\longrightarrow \mathrm{End}(H)$ and, in the PEL case,
\item a (twisted) principal polarization $\lambda: H\buildrel\sim\over\longrightarrow H^{\vee}\otimes_{\mathbb{Z}_p} \mathbb{L}$, where $H^{\vee}$ is the dual $p$-divisible group, and $\mathbb{L}$ is a $1$-dimensional smooth $\mathbb{Z}_p$-local system.
\end{altitemize}

These data are subject to the following conditions.
\begin{altenumerate}
\item[{\rm (i)}] In the PEL case, we assume that the Rosati involution $\ast_\lambda$ on $\mathrm{End}(H)$ induced by $\lambda$ is compatible with $\ast$ on $\mathcal{O}_B$, i.e $\iota(a)^{\ast_\lambda} = \iota(a^\ast)$ for all $a\in \mathcal{O}_B$.
\item[{\rm (ii)}] Locally on $S$ there is an isomorphism of $\mathcal{O}_B\otimes \mathcal{O}_S$-modules between the Lie algebra of the universal vector extension of $H$ and $\Lambda\otimes_{\mathbb{Z}_p} \mathcal{O}_S$.
\item[{\rm (iii)}] The determinant condition holds true, i.e. we have an identity of polynomial functions in $a\in \mathcal{O}_B$:
\[
\mathrm{det}_{\mathcal{O}_S}(a | \Lie X) = \mathrm{det}_E (a | V_0)\ .
\]
We refer to \cite{RapoportZinkPeriodSpaces}, 3.23, for a detailed discussion of this condition.\footnote{It is stronger than requiring just an identity of the evaluations at all $a\in \mathcal{O}_B$.}
\end{altenumerate}
\end{definition}

\begin{rem} One may wonder about the appearance of $\mathbb{L}$. First, we remark that its role is of minor importance: $\mathbb{L}$ deforms uniquely, and if one deforms an object over an algebraically closed field, then one can ignore $\mathbb{L}$ completely. This happens in particular in \cite{RapoportZinkPeriodSpaces}. Later, we will be interested in deforming $p$-divisible groups with $\mathcal{D}$-structure over finite fields, and such twisted forms will appear in the global applications: The $p$-divisible group with $\mathcal{D}$-structure over the finite field will be defined by descent from a similar object over $\bar{\mathbb{F}}_p$, and the descent datum naturally introduces this twist. We also mention that it is closely related to the notion of $c$-polarization occuring in Kottwitz' paper \cite{KottwitzPoints}, cf. Proposition \ref{VirtAbVarPDiv}.
\end{rem}

Before going on, let us consider $p$-divisible groups $\overline{\underline{H}}$ with $\mathcal{D}$-structure over perfect fields $\kappa$ of characteristic $p$. We get an associated (covariant) Dieudonn\'{e} module $(M,F)$; here $M$ is a free $W(\kappa)$-module. Moreover, $M$ carries a left action of $\mathcal{O}_B$, and an $\ast$-hermitian perfect form $M\otimes_{W(\kappa)} M \rightarrow W(\kappa)$, well-defined up to a scalar. Let us assume for the moment that $M[\frac 1p]$ is isomorphic to $V\otimes_{\mathbb{Z}_p} W(\kappa)$ as a $B$-module; this is satisfied if $\overline{\underline{H}}$ admits a deformation to a mixed characteristic discrete valuation ring as a $p$-divisible group with $\mathcal{D}$-structure, by assumption (ii) in the definition. From Lemma \ref{CrucialLemma}, it follows that in this case, we can find an isomorphism between $M$ and $\Lambda\otimes_{\mathbb{Z}_p} W(\kappa)$ as $\ast$-hermitian $\mathcal{O}_B$-modules. Let us fix such an isomorphism. The Frobenius operator $F$ on $M$ takes the form $F=p\delta \sigma$ for some $\delta\in \mathbf{G}(W(\kappa)[\frac 1p])$; moreover, changing the isomorphism changes $\delta$ by a $\sigma$-conjugate under $\mathbf{G}(W(\kappa))$. The normalization of $\delta$ is chosen to match the normalization in \cite{KottwitzPoints}. We also define $b=p\delta$, considered as a $\sigma$-conjugacy class under $\mathbf{G}(W(\bar{\kappa})[\frac 1p])$ in $\mathbf{G}(W(\bar{\kappa})[\frac 1p])$.

Let us add a word about parametrizing $\mathbb{L}$. Giving $\mathbb{L}$ is equivalent to giving $\mathbb{Q}_p/\mathbb{Z}_p\otimes_{\mathbb{Z}_p} \mathbb{L}$, which is an \'{e}tale $1$-dimensional $p$-divisible group. Its Dieudonn\'{e} module can be trivialized to $W(\kappa)$, where $F$ acts as $pd$ for some $d\in W(\kappa)^\times$. Then in general for any $p$-divisible group $X$ over $\kappa$ with Dieudonn\'{e} module $(M,F)$, the Dieudonn\'{e} module of $X\otimes \mathbb{L}$ is given by $(M,dF)$.

We are ready to state the deformation problem. Let $\underline{\overline{H}}$ be a $p$-divisible group with $\mathcal{D}$-structure over a perfect field $\kappa$ of characteristic $p$, which we give the structure of an $\mathcal{O}_E$-algebra via a fixed map $\mathcal{O}_E\longrightarrow \kappa$.

\begin{definition} Let $\mathrm{Def}_{\underline{\overline{H}}}$ be the functor that associates to every artinian local $\mathcal{O}_E$-algebra $R$ with residue field $\kappa$ the set of isomorphism classes of $p$-divisible groups with $\mathcal{D}$-structure $\underline{\tilde{H}}$ over $R$ endowed with a trivialization
\[
\underline{\overline{H}}\buildrel \cong \over \longrightarrow \underline{\tilde{H}}\otimes_R k\ .
\]
\end{definition}

\begin{thm}\label{Representability} The functor $\mathrm{Def}_{\underline{\overline{H}}}$ is pro-representable by a complete noetherian local $\mathcal{O}_E$-algebra $R_{\underline{\overline{H}}}$ with residue field $\kappa$.
\end{thm}

\begin{proof} By the results of Illusie, \cite{Illusie}, Corollaire 4.8 (i), the deformation problem for the $p$-divisible group $\overline{H}$ itself is pro-representable by a complete noetherian local $\mathcal{O}_E$-algebra $R_{\overline{H}}$ with residue field $\kappa$. By rigidity of $p$-divisible groups, the existence of liftings of the extra structure defines a quotient of $R_{\overline{H}}$. Obviously, the conditions (i) and (iii) define a further quotient $R_{\underline{\overline{H}}}$. But condition (iii) implies condition (ii) by \cite{RapoportZinkPeriodSpaces}, 3.23 c), so that $R_{\underline{\overline{H}}}$ represents $\mathrm{Def}_{\underline{\overline{H}}}$.
\end{proof}

\begin{rem} At least if $\kappa$ is algebraically closed, these deformation spaces are formal completions of the deformation spaces considered by Rapoport-Zink, \cite{RapoportZinkPeriodSpaces}.
\end{rem}

Let $k^\circ$ be the complete discrete valuation ring with residue field $\kappa$ that is unramified over $\mathcal{O}_E$ (in the sense that a uniformizer of $\mathcal{O}_E$ stays a uniformizer in $k^\circ$), and let $k$ be its fraction field. We consider the generic fibre $X_{\underline{\overline{H}}}$ of $\mathrm{Spf}\ R_{\underline{\overline{H}}}$, as a rigid-analytic space over $k$.

Associated to $\overline{\mu}$, we get a homogeneous variety $\mathcal{F}$ for $\mathbf{G}$ over $E$, cf. \cite{RapoportZinkPeriodSpaces}, 1.31. Let $d$ be its dimension. Let $\mathcal{F}^{\mathrm{rig}}$ be the associated rigid-analytic variety over $E$.

\begin{thm}\label{Smoothness} There is an \'{e}tale morphism of rigid-analytic varieties over $k$ (the \emph{period morphism})
\[
\pi: X_{\underline{\overline{H}}}\rightarrow \mathcal{F}^{\mathrm{rig}}\otimes_E k\ .
\]
In particular, $X_{\underline{\overline{H}}}$ is smooth of dimension $d$.
\end{thm}

\begin{proof} The same argument as for smoothness of Rapoport-Zink spaces applies, cf. Proposition 5.17 in \cite{RapoportZinkPeriodSpaces}: The period mapping exists in this setting by Proposition 5.15 in \cite{RapoportZinkPeriodSpaces}, and it is \'{e}tale by the same arguments.
\end{proof}

Further, we have a universal $p$-adic Tate module $T_{\underline{\overline{H}}}$ over $X_{\underline{\overline{H}}}$, with $\mathcal{D}$-structure in the obvious sense. For any compact open subgroup $K\subset K_0=\mathbf{G}(\mathbb{Z}_p)$, we let $X_{\underline{\overline{H}},K}$ be the finite \'{e}tale covering of $X_{\underline{\overline{H}}}$ parametrizing level-$K$-structures, i.e. (over each connected component of $X_{\underline{\overline{H}}}$) $\pi_1$-invariant $K$-orbits $\overline{\eta}$ of isomorphisms
\[
\eta: \Lambda\buildrel\sim\over\longrightarrow T_{\underline{\overline{H}}}
\]
which are $\mathcal{O}_B$-linear and, in the PEL case, preserve the hermitian forms up to a scalar. The existence of such isomorphisms over geometric points follows from Lemma \ref{CrucialLemma}. Together with the discussion in \cite{RapoportZinkPeriodSpaces} starting in Subsection 5.32, this implies that if $K\subset K^\prime$ is normal, the covering $X_{\underline{\overline{H}},K}/X_{\underline{\overline{H}},K^\prime}$ is Galois with Galois group $K^\prime/K$.

Now Proposition 1.20 of \cite{RapoportZinkPeriodSpaces} gives the following result.

\begin{prop} Assume that $X_{\underline{\overline{H}}}\neq \emptyset$. Then $\kappa_{\mathbf{G}}(b) = \mu^\sharp$, with notation as explained below.
\end{prop}

\begin{proof} Recall that the existence of some $x\in X_{\underline{\overline{H}}}$ implies that the (covariant) Dieudonn\'{e} module of $\underline{\overline{H}}$ can be trivialized as an $\mathcal{O}_B$-module (with $\ast$-hermitian form) to $\Lambda\otimes_{\mathbb{Z}_p} W(\kappa)$. Then the Frobenius operator $F$ defines a $\sigma$-conjugacy class $b\in B(\mathbf{G})$, where the latter denotes the set of $\sigma$-conjugacy classes in $\mathbf{G}(W(\bar{\kappa})[\frac 1p])$ (the latter set is independent of the algebraically closed field $\bar{\kappa}$ of characteristic $p$). Kottwitz, cf. \cite{KottwitzStabilization}, Section 6, constructs a map
\[
\kappa_{\mathbf{G}}: B(\mathbf{G})\rightarrow X^{\ast}(Z(\hat{\mathbf{G}})^\Gamma)\ ,
\]
where $\hat{\mathbf{G}}$ is the dual group, $Z(\hat{\mathbf{G}})$ its center, and $\Gamma$ is the absolute Galois group of $\mathbb{Q}_p$.

On the other hand, choose a cocharacter $\mu: \mathbb{G}_m\rightarrow \mathbf{G}$ representing $\pi(x)\in \mathcal{F}^{\mathrm{rig}}$. Then the pair $(b,\mu)$ is admissible in the sense of \cite{RapoportZinkPeriodSpaces}, cf. 3.19 a) of \cite{RapoportZinkPeriodSpaces}. Moreover, $\mu$ defines a character of $Z(\hat{\mathbf{G}})$, and by restriction a character $\mu^\sharp$ of $Z(\hat{\mathbf{G}})^\Gamma$; we remark that $\mu^\sharp$ depends only on the conjugacy class $\overline{\mu}$.

In general, $\kappa_{\mathbf{G}}(b) - \mu^\sharp$ lies in
\[
H^1(\mathbb{Q}_p,\mathbf{G})\cong X^{\ast}(Z(\hat{\mathbf{G}})^\Gamma)_{\mathrm{tor}}\subset X^{\ast}(Z(\hat{\mathbf{G}})^\Gamma)
\]
and measures the difference between $T_{\underline{\overline{H}},x}\otimes_{\mathbb{Z}_p} \mathbb{Q}_p$ and $V$ as $B$-modules (with hermitian form up to scalar, in the PEL case), by Proposition 1.20 of \cite{RapoportZinkPeriodSpaces}. As they are isomorphic in our case as checked above, we get the proposition.
\end{proof}

We have associated to any $p$-divisible group with $\mathcal{D}$-structure $\overline{\underline{H}}$ over $\kappa$ such that $X_{\overline{\underline{H}}}\neq \emptyset$ an element $\delta \in \mathbf{G}(W(\kappa)[\frac 1p])$, well defined up to $\sigma$-conjugation by $\mathbf{G}(W(\kappa))$, such that $p\Lambda\subset p\delta \Lambda\subset \Lambda$ and $\kappa_{\mathbf{G}}(p \delta) = \mu^\sharp$. This is summarized in the following proposition.

\begin{prop}\label{ParamPDivGroupDStrPerfField} For any perfect field $\kappa$ of characteristic $p$ which is an $\mathcal{O}_E$-algebra, the association $\underline{\overline{H}}\mapsto \delta\in \mathbf{G}(W(\kappa)[\frac 1p])$ defines an injection from the set of isomorphism classes of $p$-divisible groups with $\mathcal{D}$-structure $\underline{\overline{H}}$ over $\kappa$ such that $X_{\overline{\underline{H}}}\neq \emptyset$ into the set of $\mathbf{G}(W(\kappa))$-$\sigma$-conjugacy classes in $\mathbf{G}(W(\kappa)[\frac 1p])$ with the properties $p\Lambda\subset p\delta \Lambda\subset \Lambda$ and $\kappa_{\mathbf{G}}(p \delta) = \mu^\sharp$.
\end{prop}

Conversely, let us start with an element $\delta\in \mathbf{G}(W(\kappa)[\frac 1p])$ such that $p\Lambda\subset p\delta\Lambda\subset \Lambda$ and $\kappa_{\mathbf{G}}(p\delta) = \mu^\sharp$.  We can construct the $p$-divisible group with $\mathcal{O}_B$-action $\overline{\underline{H}}_{\delta}$ over $\kappa$ whose covariant Dieudonn\'{e} module is given by $\Lambda\otimes_{\mathbb{Z}_p} W(\kappa)$ with $F$ acting by $p\delta \sigma$. Moreover, in the PEL case, the condition $\kappa_{\mathbf{G}}(p\delta) = \mu^\sharp$ implies that for any character $\chi: \mathbf{G}\rightarrow \mathbb{G}_m$, we have $\langle \mu, \chi \rangle = \ord_p \chi(p\delta)$, cf. 3.19 b) of \cite{RapoportZinkPeriodSpaces}. Applying this for the multiplier morphism $g\mapsto \chi(g) = gg^\ast$, we get $1= \langle \mu, \chi \rangle = \ord_p \chi(p\delta)$. In particular, there is some $\mathbb{L}$ parametrized by $d=p \chi(\delta)$, cf. the discussion after Definition \ref{DefPDivGroupDStr}. This gives a twisted principal polarization
\[
\lambda: \overline{\underline{H}}_{\delta}\rightarrow \overline{\underline{H}}_{\delta}^{\vee}\otimes \mathbb{L}\ .
\]
These combine to a pair, resp. a quadruple, $\underline{H}$ satisfying conditions (i) and (ii) of a $p$-divisible group with $\mathcal{D}$-structure over $\kappa$. It does not necessarily satisfy condition (iii), however. This amounts to saying that the map in the previous proposition is not necessarily surjective.

\begin{definition} We say that $\underline{\overline{H}}$ has controlled cohomology if $X_{\underline{\overline{H}},K}$ has controlled cohomology for all normal pro-$p$ open subgroups $K\subset \mathbf{G}(\mathbb{Z}_p)$, and all $\ell\neq p$.
\end{definition}

Note that this notion depends only on the base-change of $\underline{\overline{H}}$ to an algebraic closure $\bar{\kappa}$ of $\kappa$. Moreover, it is enough to check it for a cofinal system of normal pro-$p$ open subgroups $K\subset \mathbf{G}(\mathbb{Z}_p)$, by Proposition \ref{GaloisCoverCohom}.

In our previous work, \cite{ScholzeLLC}, Theorem 2.4, we had proved an algebraization result that shows that in the EL case considered there, all $\underline{\overline{H}}$ have controlled cohomology. This made strong use of Faltings's theory of group schemes with strict $\mathcal{O}$-action.

In the case of unramified PEL data, one can use a result of Wedhorn to prove the same result.

\begin{prop}\label{AlgUnram} Assume that $F$ (equivalently, $F_0$) is unramified over $\mathbb{Q}_p$. In the PEL case, also assume that $p\neq 2$. Then all $p$-divisible groups $\underline{\overline{H}}$ with $\mathcal{D}$-structure over a perfect field $\kappa$ of characteristic $p$ have controlled cohomology. Moreover, for any normal subgroup $K\subset \mathbf{G}(\mathbb{Z}_p)$, we have
\[
H^\ast(X_{\underline{\overline{H}},K}\otimes_k \hat{\bar{k}},\mathbb{Q}_\ell)^{\mathbf{G}(\mathbb{Z}_p)} = H^\ast(X_{\underline{\overline{H}}}\otimes_k \hat{\bar{k}},\mathbb{Q}_\ell) = \mathbb{Q}_\ell\ .
\]
\end{prop}

\begin{proof} In the unramified case, conditions (ii) and (iii) in the definition of a $p$-divisible group with $\mathcal{D}$-structure reduce to numerical conditions on the ranks of certain locally free sheaves and hence can be checked on geometric points and are irrelevant for deformation problems. Now Theorem 2.8 and 2.15 of \cite{WedhornOortStrata} show that $R_{\overline{\underline{H}}}$ is formally smooth and the prorepresentable hull of a finitely presented functor (given as the deformation functor of any truncation $\underline{\overline{H}}[p^m]$). Now using Artin's algebraization theorem as in \cite{ScholzeLLC}, proof of Theorem 2.4, one constructs an algebraization as in Theorem \ref{AlgebraizationControlled}. Theorem \ref{AlgebraizationControlled} now shows that $\underline{\overline{H}}$ has controlled cohomology, and also gives the desired description of the unramified part of the cohomology.
\end{proof}

\begin{prop}\label{FactEL1} Assume given data of EL type, and assume that $F$ factors as a product of fields $F=\prod_i F_i$. Accordingly, all other data split into a product, which we indicate by writing $\mathcal{D} = \prod_i \mathcal{D}_i$. Similarly, a $p$-divisible group with $\mathcal{D}$-structure $\underline{H}$ over $S$ factors into a product $\underline{H} = \prod_i \underline{H}_i$, where $\underline{H}_i$ is a $p$-divisible group with $\mathcal{D}_i$-structure over $S$. The reflex field $E\subset \bar{\mathbb{Q}}_p$ is the compositum of the reflex fields $E_i\subset \bar{\mathbb{Q}}_p$.

Let $\kappa$ be a perfect field as above, and let $\underline{\overline{H}} = \prod_i \underline{\overline{H}}_i$ be a $p$-divisible group with $\mathcal{D}$-structure over $S$, given as a product of $p$-divisible groups with $\mathcal{D}_i$-structure $\underline{\overline{H}}_i$ over $S$. Accordingly,
\[
\delta=\prod_i \delta_i\in \mathbf{G}(W(\kappa)[\frac 1p])=\prod_i \mathbf{G}_i(W(\kappa)[\frac 1p])\ .
\]

Moreover, if $K=\prod_i K_i\subset \mathbf{G}(\mathbb{Z}_p)$, there are product decompositions $X_{\underline{\overline{H}},K} = \prod_i X_{\underline{\overline{H}}_i,K_i}$, the product being taken over $\mathrm{Sp}\ k$, with $k$ as above.

In particular, $\underline{\overline{H}}$ has controlled cohomology if and only if all $\underline{\overline{H}}_i$ have controlled cohomology.
\end{prop}

\begin{proof} Easy and left to reader; use Proposition \ref{ProductControlled} for the last assertion.
\end{proof}

\begin{prop}\label{ReductionToEL1} Assume given data $\mathcal{D}$ of quasi-EL type, with corresponding data $\mathcal{D}_0$ of EL type.

Then giving a $p$-divisible group with $\mathcal{D}$-structure $\underline{H} = (H,\iota,\lambda,\mathbb{L})$ over $S$ is equivalent to giving a $p$-divisible group with $\mathcal{D}_0$-structure $\underline{H}_0$ and a $1$-dimensional $\mathbb{Z}_p$-local system $\mathbb{L}$ over $S$, the correspondence being given by $H=H_0\times H_0^\vee\otimes \mathbb{L}$ with the tautological $\mathcal{O}_B = \mathcal{O}_{B_0}\times \mathcal{O}_{B_0}^\ast$-action, and the tautological twisted principal polarization
\[
\lambda: H=H_0\times H_0^\vee\otimes \mathbb{L}\rightarrow H^{\vee}\otimes \mathbb{L}= H_0^\vee\otimes \mathbb{L}\times H_0\ .
\]

In particular, let $\kappa$ and $\underline{\overline{H}}$ be as above, with associated $\underline{\overline{H}}_0$ and $\mathbb{L}$. Accordingly,
\[
\delta = (\delta_0,p^{-1} d)\in \mathbf{G}(W(\kappa)[\frac 1p]) = \mathbf{G}_0(W(\kappa)[\frac 1p])\times W(\kappa)[\frac 1p]^\times\ ,
\]
where $\mathbb{L}$ corresponds to $d$ as above. Then for
\[
K=K_0\times (1+p^m \mathbb{Z}_p)\subset \mathbf{G}(\mathbb{Z}_p) = \mathbf{G}_0(\mathbb{Z}_p)\times \mathbb{Z}_p^\times\ ,
\]
with $m\geq 1$, there is a product decomposition
\[
X_{\underline{\overline{H}},K} = X_{\underline{\overline{H}}_0,K_0}\times_{\mathrm{Sp}\, k} X_{\mathbb{L},m}\ ,
\]
where $X_{\mathbb{L},m}$ parametrizes isomorphisms between $\mathbb{L}\otimes \mu_{p^m}$ and $\mathbb{Z}/p^m\mathbb{Z}$.

In particular, $\overline{\underline{H}}$ has controlled cohomology if and only if $\underline{\overline{H}}_0$ has controlled cohomology.
\end{prop}

\begin{proof} Easy and left to reader; note that $X_{\mathbb{L},m}$ is finite over $\mathrm{Sp}\, k$ and hence has controlled cohomology.
\end{proof}

Finally, we will need a continuity statement about the action of the automorphism group of $\underline{\overline{H}}$ on the cohomology of $X_{\underline{\overline{H}},K}$ in the case that $\underline{\overline{H}}$ has controlled cohomology.

\begin{prop}\label{JActionContinuous} Assume that $\underline{\overline{H}}$ has controlled cohomology. Then for any normal pro-$p$ open subgroup $K\subset \mathbf{G}(\mathbb{Z}_p)$, there is an integer $m\geq 1$ such that for all automorphisms $j$ of $\underline{\overline{H}}$ that act trivially on $\overline{H}[p^m]$, the induced action on
\[
H^i(X_{\underline{\overline{H}},K}\otimes_k \hat{\bar{k}},\mathbb{Q}_\ell)
\]
is trivial for all $i$.
\end{prop}

\begin{proof} It suffices to check this for a cofinal system of $K$. We take $K$ as the kernel of the projection $\mathbf{G}(\mathbb{Z}_p)\rightarrow \mathbf{G}(\mathbb{Z}/p^{m_1}\mathbb{Z})$. In this case $X_{\underline{\overline{H}},K}/X_{\underline{\overline{H}}}$ parametrizes $\mathcal{O}_B$-linear isomorphisms $\Lambda/p^{m_1}\cong T_{\underline{\overline{H}}}/p^{m_1}$ that preserve the hermitian form up to a scalar.

Fix generators $\lambda_1,...,\lambda_{m_2}$ of $\Lambda$ as an $\mathcal{O}_B$-module. Let $\underline{H} = (H,\ldots)$ be the universal deformation of $\underline{\overline{H}}$ over $R_{\underline{\overline{H}}}$. Consider the finite flat cover $H[p^{m_1}] / \mathrm{Spf}\ R_{\underline{\overline{H}}}$, and let $H[p^{m_1}]_\eta / X_{\underline{\overline{H}}}$ be its generic fibre as a rigid-analytic variety over $k$. Let $Y = H[p^{m_1}]_\eta^{m_2/X_{\underline{\overline{H}}}}$ be the $m_2$-fold fibre product over $X_{\underline{\overline{H}}}$. Note that $Y$ is the generic fibre of a complete noetherian semilocal $k^\circ$-algebra $R$ whose residue fields are finite extensions of $\kappa$: Let $H[p^{m_1}] = \mathrm{Spf}\ R_1$; then
\[
R = R_1\otimes_{R_{\underline{\overline{H}}}} \cdots \otimes_{R_{\underline{\overline{H}}}} R_1\ ,
\]
which is finite flat (in fact free) over $R_{\underline{\overline{H}}}$. We have a closed immersion
\[
X_{\underline{\overline{H}},K}\rightarrow Y = H[p^{m_1}]_\eta^{m_2/X_{\underline{\overline{H}}}}
\]
by sending the isomorphism $\Lambda/p^{m_1} \cong T/p^{m_1}$ to the images of $\lambda_1,...,\lambda_{m_2}$.

Because $\underline{\overline{H}}$ has controlled cohomology, there is some quasicompact admissible open subset $X_{\underline{\overline{H}},K,m_3}\subset X_{\underline{\overline{H}},K}$ with the same cohomology, where we recall that the rigid-analytic variety $X_{\underline{\overline{H}},K}$ is defined as the union of quasicompact admissible open subsets $X_{\underline{\overline{H}},K,n}$. We use Proposition \ref{Continuity} to produce an admissible open neighborhood $U\subset X_{\underline{\overline{H}},K,m_3}\times X_{\underline{\overline{H}},K,m_3}$ of the diagonal such that any automorphism of $X_{\underline{\overline{H}},K,m_3}$ whose graph is contained in $U$ acts trivially on the cohomology.

Because $X_{\underline{\overline{H}},K,m_3}$ is quasicompact, there is some $m_4$ such that $X_{\underline{\overline{H}},K,m_3}\subset Y_{m_4}$, where $Y=\bigcup_m Y_m$ as usual. Moreover, there is an admissible open neighborhood $V\subset Y_{m_4}\times Y_{m_4}$ of the diagonal such that
\[
V\cap \left(X_{\underline{\overline{H}},K,m_3}\times X_{\underline{\overline{H}},K,m_3}\right)\subset U\ .
\]

Summarizing, it suffices to find an integer $m\geq 1$ such that for any automorphism $j$ of $\underline{\overline{H}}$ that acts trivially on $\overline{H}[p^m]$, the graph $\Gamma_j\subset Y_{m_4}\times Y_{m_4}$ of $j$ acting on $Y_{m_4}$ is contained in $V$.

Using Lemma \ref{ContinuityHelp}, it suffices to prove that for any open ideal $I\subset R$, there is some integer $m\geq 1$ such that any automorphism $j$ of $\underline{\overline{H}}$ that acts trivially on $\overline{H}[p^m]$, also acts trivially on $R/I$. This reduces us to proving the statement with $R_1$ in place of $R$.

Let $\mathfrak{m}_{\underline{\overline{H}}}\subset R_{\underline{\overline{H}}}$ be the maximal ideal. It suffices to check that for any $n\geq 1$ there is some $m\geq 1$ such that any automorphism $j$ of $\underline{\overline{H}}$ acting trivially on $\overline{H}[p^m]$ lifts to an automorphism of $\underline{H}\otimes_{R_{\underline{\overline{H}}}} R_{\underline{\overline{H}}}/ \mathfrak{m}_{\underline{\overline{H}}}^n$ that is trivial on $p^{m_1}$-torsion.

By induction, it suffices to prove that there is some integer $m_1^\prime\geq m_1$ such that any automorphism $j$ of $\underline{\overline{H}}$ that lifts to an automorphism of $\underline{H}\otimes_{R_{\underline{\overline{H}}}} R_{\underline{\overline{H}}}/ \mathfrak{m}_{\underline{\overline{H}}}^{n-1}$ trivial on $p^{m_1^\prime}$-torsion, lifts further to an automorphism of $\underline{H}\otimes_{R_{\underline{\overline{H}}}} R_{\underline{\overline{H}}}/ \mathfrak{m}_{\underline{\overline{H}}}^n$ trivial on $p^{m_1}$-torsion.

The assertion is equivalent to the existence of lifts of $\frac{j-1}{p^{m_1}}, \frac{j^{-1}-1}{p^{m_1}}$ acting on $\underline{H}\otimes_{R_{\underline{\overline{H}}}} R_{\underline{\overline{H}}}/ \mathfrak{m}_{\underline{\overline{H}}}^{n-1}$. By \cite{Illusie}, Corollaire 4.3 b), for any $p$-divisible groups $G_1$, $G_2$ over $R_{\underline{\overline{H}}}/ \mathfrak{m}_{\underline{\overline{H}}}^n$ with restriction $G_1^\prime$, $G_2^\prime$ to $R_{\underline{\overline{H}}}/ \mathfrak{m}_{\underline{\overline{H}}}^{n-1}$, the existence of a lift of some homomorphism $f: G_1^\prime\rightarrow G_2^\prime$ to a homomorphism $G_1\rightarrow G_2$ is equivalent to the existence of a lift of $f[p]: G_1^\prime[p]\rightarrow G_2^\prime[p]$ to a homomorphism $G_1[p]\rightarrow G_2[p]$. But taking $m_1^\prime=m_1+1$, we know that $f[p]$ is just the zero morphism for the two homomorphisms $f=\frac{j-1}{p^{m_1}}, \frac{j^{-1}-1}{p^{m_1}}$ we are interested in. These obviously lift, finishing the proof.
\end{proof}

\section{Definition of the test function}\label{DefTestFunction}

Let $I_E\subset W_E$ be the inertia and Weil group of $E$, and fix a geometric Frobenius element $\mathrm{Frob}\in W_E$. Fix some integer $j\geq 1$. Our aim is to define a function $\phi_{\tau,h}\in C_c^{\infty}(\mathbf{G}(\mathbb{Q}_{p^r}))$ depending on an element $\tau\in \mathrm{Frob}^j I_E\subset W_E$ and a function $h\in C_c^{\infty}(\mathbf{G}(\mathbb{Z}_p))$ with values in $\mathbb{Q}$. Here we set $r=j[\kappa_E:\mathbb{F}_p]$, where $\kappa_E$ is the residue field of $E$.

We regard $\mathbb{F}_{p^r}$ as the degree-$j$-extension of $\kappa_E$; in particular, it is an $\mathcal{O}_E$-algebra. Fix the Haar measures on $\mathbf{G}(\mathbb{Q}_p)$, resp. $\mathbf{G}(\mathbb{Q}_{p^r})$, that give $\mathbf{G}(\mathbb{Z}_p)$, resp. $\mathbf{G}(\mathbb{Z}_{p^r})$, volume $1$.

\begin{definition} Let $\delta\in \mathbf{G}(\mathbb{Q}_{p^r})$. Define
\[
\phi_{\tau,h}(\delta) = 0
\]
unless $\delta$ is associated to some $p$-divisible group with $\mathcal{D}$-structure $\underline{\overline{H}}$ over $\mathbb{F}_{p^r}$ under the correspondence of Proposition \ref{ParamPDivGroupDStrPerfField}. In the latter case, assume first that $\underline{\overline{H}}$ has controlled cohomology. Then define
\[
\phi_{\tau,h}(\delta) = \mathrm{tr}(\tau\times h|H^{\ast}(X_{\underline{\overline{H}},K}\otimes_k \hat{\bar{k}}, \mathbb{Q}_\ell))\ ,
\]
for any normal compact pro-$p$ open subgroup $K\subset K_0$ such that $h$ is $K$-biinvariant. If $\underline{\overline{H}}$ does not have controlled cohomology, define $\phi_{\tau,h}(\delta)=0$.
\end{definition}

\begin{prop} The function $\phi_{\tau,h}: \mathbf{G}(\mathbb{Q}_{p^r})\rightarrow \mathbb{Q}_{\ell}$ is well-defined and takes values in $\mathbb{Q}$ independent of $\ell$. Its support is contained in the compact set of all $\delta\in \mathbf{G}(\mathbb{Q}_{p^r})$ satisfying $p\Lambda\subset p\delta\Lambda\subset\Lambda$ and $\kappa_{\mathbf{G}}(p\delta) = \mu^\sharp$.
\end{prop}

\begin{proof} The last point is clear. We have to see that any two choices $K_1\subset K_2$ of $K$ give the same value. This follows from Proposition \ref{GaloisCoverCohom} as follows:
\[\begin{aligned}
&\phantom{=} \mathrm{tr}(\tau\times h|H^{\ast}(X_{\underline{\overline{H}},K_2}\otimes_k \hat{\bar{k}}, \mathbb{Q}_\ell))\\
&= \mathrm{vol}(K_2)^{-1} \sum_{g\in \mathbf{G}(\mathbb{Z}_p)/K_2} h(g) \mathrm{tr}(\tau\times g|H^{\ast}(X_{\underline{\overline{H}},K_2}\otimes_k \hat{\bar{k}}, \mathbb{Q}_\ell))\\
&=  \mathrm{vol}(K_2)^{-1} \sum_{g\in \mathbf{G}(\mathbb{Z}_p)/K_2} h(g) \mathrm{tr}(\tau\times g|H^{\ast}(X_{\underline{\overline{H}},K_1}\otimes_k \hat{\bar{k}}, \mathbb{Q}_\ell)^{K_2/K_1})\\
&=  \mathrm{vol}(K_1)^{-1} \sum_{g\in \mathbf{G}(\mathbb{Z}_p)/K_1} h(g) \mathrm{tr}(\tau\times g|H^{\ast}(X_{\underline{\overline{H}},K_1}\otimes_k \hat{\bar{k}}, \mathbb{Q}_\ell))\\
&= \mathrm{tr}(\tau\times h|H^{\ast}(X_{\underline{\overline{H}},K_1}\otimes_k \hat{\bar{k}}, \mathbb{Q}_\ell))\ .
\end{aligned}\]

To get the independence of $\ell$, we note that it suffices to prove the independence of $\ell$ of
\[
\mathrm{tr}(\tau\times g|H^{\ast}(X_{\underline{\overline{H}},K}\otimes_k \hat{\bar{k}}, \mathbb{Q}_\ell))
\]
for all $K$ and $g\in \mathbf{G}(\mathbb{Z}_p)$, assuming $\underline{\overline{H}}$ has controlled cohomology. In particular, we assume that $X_{\underline{\overline{H}},K}$ has controlled cohomology, so we may replace it by some $\mathbf{G}(\mathbb{Z}_p)$-invariant quasicompact open subset $U=X_{\underline{\overline{H}},K,m}\subset X_{\underline{\overline{H}},K}$ with the same cohomology. Moreover, we can twist $U$ by the unramified action of the absolute Galois group of $k$ sending a geometric Frobenius element to $g$, which acts via the finite quotient $\mathbf{G}(\mathbb{Z}_p)/K$; this gives a quasicompact smooth separated rigid variety $V$ over $k$ such that $U\otimes_k \hat{\bar{k}}\cong V\otimes_k \hat{\bar{k}}$, and with the action of $\tau\times g$ on the left-hand side corresponding to the action of $\tau$ on the right-hand side. In particular,
\[
\mathrm{tr}(\tau\times g|H^{\ast}(X_{\underline{\overline{H}},K}\otimes_k \hat{\bar{k}}, \mathbb{Q}_\ell)) = \mathrm{tr}(\tau | H^\ast(V\otimes_k \hat{\bar{k}},\mathbb{Q}_{\ell}))\ .
\]
The latter term is independent of $\ell$ by Theorem 7.1.10 of \cite{MiedaIndependence}.
\end{proof}

\begin{prop} The function $\phi_{\tau,h}$ is locally constant, so that it defines an element $\phi_{\tau,h}\in C_c^{\infty}(\mathbf{G}(\mathbb{Q}_{p^r}))$.
\end{prop}

\begin{proof} Take any element $\delta\in \mathbf{G}(\mathbb{Q}_{p^r})$; we want to find a small open neighborhood $U$ of $\delta$ such that $\phi_{\tau,h}(\delta^\prime) = \phi_{\tau,h}(\delta)$ for all $\delta^\prime\in U$. The conditions $\kappa_{\mathbf{G}}(p\delta) = \mu^\sharp$ and $p\Lambda\subset p\delta \Lambda\subset \Lambda$ define an open and closed subset outside of which the function vanishes identically. In particular, we may assume that $\delta$ satisfies these conditions.

The construction after Proposition \ref{ParamPDivGroupDStrPerfField} constructs a pair $\overline{\underline{H}}=(\overline{H},\iota)$ (resp. quadruple $\overline{\underline{H}}=(\overline{H},\iota,\lambda,\mathbb{L})$ in the PEL case), which satisfies all conditions of being a $p$-divisible group with $\mathcal{D}$-structure over $\mathbb{F}_{p^r}$ except possibly the determinant condition, i.e. condition (iii) of Definition \ref{DefPDivGroupDStr}.

We want to see that if $\delta^\prime$ is sufficiently close to $\delta$, then over $\bar{\mathbb{F}}_p$, the associated $\overline{\underline{H}}_{\delta}$ and $\overline{\underline{H}}_{\delta^\prime}$ become isomorphic. This follows from the following lemma.

\begin{lem} Let $\kappa$ be an algebraically closed field of characteristic $p$, let $L=W(\kappa)[\frac 1p]$, and let $G$ be any linear algebraic group over $L$. Then for any $b\in G(L)$, the map $G(L)\rightarrow G(L)$ mapping $g$ to $g^{-1} bg^\sigma$ is open.
\end{lem}

\begin{proof} By standard arguments, it is enough to check the statement on the Lie algebra, which reads: The map $\mathfrak{g}\rightarrow \mathfrak{g}$ mapping $x$ to $-x + (\mathrm{Ad} b)(x^\sigma)$ is open, where $\mathfrak{g}$ denotes the Lie algebra of $G$. Identifying $\mathfrak{g}$ with $L^n$ for some $n$, this follows from the following lemma.

\begin{lem} Let $A\in \mathrm{GL}_n(L)$. Then the map $L^n\rightarrow L^n$ mapping $x$ to $x - Ax^\sigma$ is open.
\end{lem}

\begin{rem} By $\mathbb{Q}_p$-linearity, it follows that the map is also surjective.
\end{rem}

\begin{proof} It is easily seen that one may replace $A$ by a $\sigma$-conjugate. Using the Dieudonn\'{e}-Manin classification, and reducing to simple factors, we may thus assume that $A$ has the form
\[
A = \left(\begin{array}{ccccc} 0 & 1 & 0 & \cdots & 0 \\ 0 & 0 & 1 & \cdots & 0 \\ \vdots & \vdots & \ddots  & \ddots & \vdots \\ 0 & 0 & 0 & \cdots & 1 \\ p^k & 0 & 0 & \cdots & 0 \end{array}\right)\ ,
\]
for some integer $k$ prime to the size of the matrix. We distinguish the cases $k>0$, $k=0$ and $k<0$. If $k>0$, then $A$ is topologically nilpotent and we can give the inverse of $x-Ax^\sigma$ explicitly as $x+Ax^\sigma + A^2x^{\sigma^2} + \ldots$. Similarly, if $k<0$, then $A^{-1}$ is topologically nilpotent, and we can give the inverse of $x-Ax^\sigma$ explicitly as $-A^{-1}x^{\sigma^{-1}} - A^{-2}x^{\sigma^{-2}} - \ldots$. Finally, if $k=0$, then $A=1$, and it suffices to show that the $\mathbb{Z}_p$-linear map $W(\kappa)\rightarrow W(\kappa)$, $x\mapsto x-x^\sigma$, is surjective. This can be checked modulo $p$, where it follows directly from the fact that $\kappa$ is algebraically closed.
\end{proof}
\end{proof}

We use this lemma for $\delta\in \mathbf{G}(\mathbb{Q}_{p^r})\subset \mathbf{G}(W(\bar{\mathbb{F}}_p)[\frac 1p])$. Let $m\geq 1$ be some positive integer (to be chosen later). Let $V$ be the open neighborhood of $\delta$ given as the image of $\mathrm{ker}(\mathbf{G}(W(\bar{\mathbb{F}}_p))\rightarrow \mathbf{G}(W(\bar{\mathbb{F}}_p)/p^m)))$ under $g\mapsto g^{-1} \delta g^\sigma$, and let $U=V\cap \mathbf{G}(\mathbb{Q}_{p^r})$.

It follows that for any $\delta^\prime\in U$, the $p$-divisible groups $\overline{\underline{H}}_{\delta}$ and $\overline{\underline{H}}_{\delta^\prime}$ associated to $\delta$, resp. $\delta^\prime$ become isomorphic over $\bar{\mathbb{F}}_p$. In particular, one satisfies the determinant condition if and only if the other one does, and in this case we get isomorphisms
\[
X_{\overline{\underline{H}}_{\delta},K}\otimes_k \hat{\bar{k}}\cong X_{\overline{\underline{H}}_{\delta^\prime},K}\otimes_k \hat{\bar{k}}\ ,
\]
for all open subgroups $K\subset \mathbf{G}(\mathbb{Z}_p)$, so that $\underline{\overline{H}}_{\delta}$ has controlled cohomology if and only if $\underline{\overline{H}}_{\delta^\prime}$ has controlled cohomology.

In particular, we get the desired statement that $\phi_{\tau,h}$ is locally constant at $\delta$ unless $\overline{\underline{H}}_\delta$ is a $p$-divisible group with $\mathcal{D}$-structure and has controlled cohomology. In this final case, we choose $m$ large enough such that Proposition \ref{JActionContinuous} applies. Twisting with elements of $\mathrm{ker}(\mathbf{G}(W(\bar{\mathbb{F}}_p))\rightarrow \mathbf{G}(W(\bar{\mathbb{F}}_p)/p^m)))$ does not change $\overline{\underline{H}}_\delta[p^m]$, so that the (different) actions of $\tau$ on
\[
X_{\overline{\underline{H}}_{\delta},K}\otimes_k \hat{\bar{k}}\cong X_{\overline{\underline{H}}_{\delta^\prime},K}\otimes_k \hat{\bar{k}}
\]
differ by the action of some $j$ in the automorphism group of $\overline{\underline{H}}_{\delta}$ over $\bar{\mathbb{F}}_p$ which is trivial on $p^m$-torsion points. Now Proposition \ref{JActionContinuous} implies that both actions of $\tau$ become identical on the cohomology. We find that indeed $\phi_{\tau,h}$ is locally constant.
\end{proof}

We end this section by stating some easy lemmas about $\phi_{\tau,h}$.

\begin{prop} Assume that $F$ is unramified, and $p\neq 2$ in the PEL case. Let $h$ be the idempotent associated to $\mathbf{G}(\mathbb{Z}_p)$. Then $\phi_{\tau,h}$ is the characteristic function of the double coset $\mathbf{G}(\mathbb{Z}_{p^r})p^{-1}\sigma(\mu(p))\mathbf{G}(\mathbb{Z}_{p^r})$, where $\mu: \mathbb{G}_m\rightarrow \mathbf{G}_{\mathbb{Q}_{p^r}}$ is some representative of $\overline{\mu}$ that factors over $T_{\mathbb{Q}_{p^r}}$ for some maximal $\mathbb{Z}_{p^r}$-split torus $T\subset \mathbf{G}_{\mathbb{Z}_{p^r}}$.
\end{prop}

\begin{rem} Such $\mu$ exist by \cite{KottwitzTO}, Lemma 1.1.3 (a). The proposition implies that in the unramified case, our function $\phi_{\tau,h}(\delta)$ is equal to $\phi_r(\sigma^{-1}(\delta))$, where $\phi_r$ is the function used by Kottwitz in \cite{KottwitzPoints}. Their twisted orbital integrals obviously agree, as $\delta$ is $\sigma$-conjugate to $\sigma^{-1}(\delta)$.
\end{rem}

\begin{proof} From Proposition \ref{AlgUnram}, we know that all $p$-divisible groups with $\mathcal{D}$-structure have controlled cohomology, and that in this case
\[
\phi_{\tau,h}(\delta) = 1\ .
\]
It remains to classify the set of $\delta$ that give rise to a $p$-divisible group with $\mathcal{D}$-structure. This is done by Kottwitz, \cite{KottwitzPoints}, pages 430 -- 431.
\end{proof}

\begin{prop}\label{FactEL2} In the situation of Proposition \ref{FactEL1}, let $h$ be of the form $h=\prod h_i$ with $h_i\in C_c^{\infty}(\mathbf{G}_i(\mathbb{Z}_p))$. Then for all $\delta = \prod \delta_i\in \mathbf{G}(\mathbb{Q}_{p^r})=\prod \mathbf{G}_i(\mathbb{Q}_{p^r})$,
\[
\phi_{\tau,h}(\delta) = \prod \phi_{\tau,h_i}(\delta_i)\ .
\]
\end{prop}

\begin{proof} Follows directly from Proposition \ref{FactEL1}.
\end{proof}

\begin{prop}\label{ReductionToEL2} In the situation of Proposition \ref{ReductionToEL1}, let $h$ be of the form $h_0\times h_{\mathbb{G}_m}$, where $h_0\in C_c^\infty(\mathbf{G}_0(\mathbb{Z}_p))$ and $h_{\mathbb{G}_m}\in C_c^\infty(\mathbb{Z}_p^\times)$. Then for all
\[
\delta=(\delta_0,\delta_{\mathbb{G}_m})\in \mathbf{G}(\mathbb{Q}_{p^r}) = \mathbf{G}_0(\mathbb{Q}_{p^r})\times \mathbb{Q}_{p^r}^\times\ ,
\]
there is a factorization
\[
\phi_{\tau,h}(\delta) = \phi_{\tau,h_0}(\delta_0) \phi_{\tau,h_{\mathbb{G}_m}}(\delta_{\mathbb{G}_m})\ ,
\]
where $\phi_{\tau,h_{\mathbb{G}_m}}$ is the function with support on $p^{-1}\mathbb{Z}_{p^r}^\times$ defined by
\[
\phi_{\tau,h_{\mathbb{G}_m}}(\delta_{\mathbb{G}_m}) = h(\mathrm{Art}_{\mathbb{Q}_p}(\tau) N\delta_{\mathbb{G}_m})\ ,
\]
where $\mathrm{Art}_{\mathbb{Q}_p}: W_{\mathbb{Q}_p}\rightarrow \mathbb{Q}_p^\times$ is the local reciprocity map sending a geometric Frobenius element to a uniformizer.
\end{prop}

\begin{rem} We remark that $\phi_{\tau,h_{\mathbb{G}_m}}$ satisfies the following spectral identity. For all characters $\chi: \mathbb{Q}_p^\times\rightarrow \mathbb{C}^\times$,
\[
\tr(\phi_{\tau,h_{\mathbb{G}_m}}|\chi\circ \mathrm{Norm}_{\mathbb{Q}_{p^r}/\mathbb{Q}_p}) = \tr(\tau^{-1}|\chi \circ \mathrm{Art}_{\mathbb{Q}_p})\tr(h|\chi)\ .
\]
\end{rem}

\begin{proof} This follows from Proposition \ref{ReductionToEL1}, once one checks that
\[
\mathrm{tr}(\tau\times h_{\mathbb{G}_m} | H^\ast(X_{\mathbb{L},m}\otimes_k \hat{\bar{k}},\mathbb{Q}_\ell) ) = \phi_{\tau,h_{\mathbb{G}_m}}(\delta_{\mathbb{G}_m})
\]
for all $m$ large enough.
\end{proof}

\section{PEL Shimura varieties}\label{PELShimuraVarieties}

We first recall the definition of the PEL Shimura varieties that are also considered in Kottwitz' article \cite{KottwitzPoints}. They are associated to the following data.\\

{\bf Global PEL data.} Let $B$ be a simple $\mathbb{Q}$-algebra with center $F$ and maximal $\mathbb{Z}_{(p)}$-order $\mathcal{O}_B$ that is stable under a positive involution $\ast$ on $B$. Let $V$ be a finitely generated left $B$-module with a nondegenerate $\ast$-hermitian form $(\, ,\, )$. We assume that there is an $\mathcal{O}_B$-stable selfdual $\mathbb{Z}_{(p)}$-lattice $\Lambda\subset V$, which we fix. Moreover, we let $F_0 = F^{\ast = 1}$, which is a totally real field. We assume that at all places above $p$, the extension $F/F_0$ is unramified and the $F$-algebra $B$ is split.

We let $C=\mathrm{End}_B(V)$, and $\mathcal{O}_C = \mathrm{End}_{\mathcal{O}_B}(\Lambda)$; both carry an involution $\ast$ induced from $(\, ,\, )$. We recall, cf. \cite{KottwitzPoints}, p. 375, that over $\bar{\mathbb{Q}}$, the algebra $C$ together with the involution $\ast$ is of one of the following types:
\begin{altenumerate}
\item[(A)] $M_n\times M_n^{\mathrm{opp}}$ with $(x,y)^\ast = (y,x)$,
\item[(C)] $M_{2n}$ with $x^\ast$ being the adjoint of $x$ with respect to a nondegenerate alternating form in $2n$ variables,
\item[(D)] $M_{2n}$ with $x^\ast$ being the adjoint of $x$ with respect to a nondegenerate symmetric form in $2n$ variables.
\end{altenumerate}
We assume that case $A$ or $C$ occurs. We get the reductive group $\mathbf{G}/\mathbb{Q}$ of $B$-linear similitudes of $V$; in fact, we can extend it to an algebraic group over $\mathbb{Z}_{(p)}$ as the group representing the functor
\[
\mathbf{G}(R) = \{g\in (\mathcal{O}_C\otimes_{\mathbb{Z}_p} R)^\times \mid gg^\ast\in R^\times\}\ .
\]
Finally, we fix a homomorphism $\mathbf{h}_0: \mathbb{C}\rightarrow C\otimes \mathbb{R}$ such that $\mathbf{h}_0(\overline{z})=\mathbf{h}_0(z)^\ast$ for all $z\in \mathbb{C}$, and such that the symmetric real-valued bilinear form $(v,\mathbf{h}_0(i)w)$ on $V\otimes \mathbb{R}$ is positive definite.
\\

We write $\mathbf{h}$ for the map $\mathbb{S}\rightarrow \mathbf{G}\otimes \mathbb{R}$ from Deligne's torus $\mathbb{S}$ (i.e., the algebraic torus over $\mathbb{R}$ with $\mathbb{S}(\mathbb{R}) = \mathbb{C}^\times$) that is given on $\mathbb{R}$-valued points by $\mathbf{h}(z) = \mathbf{h}_0(z)$, $z\in\mathbb{C}^\times$. Then one gets a tower $\mathrm{Sh}_K$, $K\subset \mathbf{G}(\mathbb{A}_f)$ running through compact open subgroups of the finite adelic points of $\mathbf{G}$, of Shimura varieties associated to the pair $(\mathbf{G},\mathbf{h}^{-1})$.

A priori, these Shimura varieties are defined over $\mathbb{C}$, but they have canonical models over their reflex field $E\subset \mathbb{C}$, a finite extension of $\mathbb{Q}$. We recall that $E$ is the field of definition of the conjugacy class of the cocharacter $\mu = \mu_{\mathbf{h}}: \mathbb{G}_m\rightarrow \mathbf{G}\otimes \mathbb{C}$ given as the product of the central morphism $\mathbb{G}_m\rightarrow \mathbf{G}$ sending $t\in \mathbb{G}_m$ to multiplication by $t$ on $V$, and the composite morphism
\[
\mathbb{G}_m\rightarrow \mathbb{S}\otimes \mathbb{C}\buildrel {\mathbf{h}^{-1}\otimes \mathbb{C}}\over\longrightarrow \mathbf{G}\otimes \mathbb{C}\ .
\]
Recall that $\mathbb{S}\otimes \mathbb{C} = \mathbb{G}_m\times \mathbb{G}_m$, the first factor corresponding to the identity morphism $\mathbb{C}\rightarrow \mathbb{C}$, the second to complex conjugation $\mathbb{C}\rightarrow \mathbb{C}$. The first morphism $\mathbb{G}_m\rightarrow \mathbb{S}\otimes \mathbb{C}$ is the one coming from the identity morphism.

We note that this is not the same $\mu$ as the $\mu_{\mathrm{Kw}}$ considered by Kottwitz, which is given as the composition
\[
\mathbb{G}_m\rightarrow \mathbb{S}\otimes \mathbb{C}\buildrel {\mathbf{h}\otimes \mathbb{C}}\over\longrightarrow \mathbf{G}\otimes \mathbb{C}\ .
\]
This means that the product $\mu \mu_{Kw}$ is the central morphism $\mathbb{G}_m\rightarrow \mathbf{G}$. We choose this alternative normalization as it is the one compatible with our local normalization, which in turn is the one used by Rapoport and Zink in their book \cite{RapoportZinkPeriodSpaces}. Fix a prime $\mathfrak{p}$ of $E$ above $p$, and let $\mathcal{O}_{E_{\mathfrak{p}}}$ be the complete local ring at $\mathfrak{p}$. This allows to regard $\mu$ as a conjugacy class of cocharacters $\mu: \mathbb{G}_m\rightarrow \mathbf{G}_{\bar{\mathbb{Q}}_p}$, where $\bar{\mathbb{Q}}_p$ is an algebraic closure of $E_{\mathfrak{p}}$. In particular, after base-changing everything to $\mathbb{Q}_p$, we get data $\mathcal{D}$ of PEL type as defined in Section \ref{DefSpaces}.

Next, we recall that (finite disjoint unions of) these Shimura varieties can be described as moduli spaces of abelian varieties with polarization, endomorphism, and level structure. This also leads to integral models of these Shimura varieties.

\begin{definition} Let $K^p\subset \mathbf{G}(\mathbb{A}_f^p)$ be a sufficiently small compact open subgroup. Let $\mathfrak{M}_{K^p}$ be the contravariant set-valued functor on the category of locally noetherian schemes $S$ over $\mathcal{O}_{E_{\mathfrak{p}}}$ associating to $S$ the set of isomorphism classes of quadruples $(A,\iota,\lambda,\overline{\eta})$, consisting of
\begin{altitemize}
\item an abelian scheme $A$ up to prime-to-$p$-isogeny over $S$,
\item an action $\iota: \mathcal{O}_B\rightarrow \mathrm{End}(A)$,
\item a prime-to-$p$-isogeny $\lambda: A\rightarrow A^\vee$ which is a polarization,
\item a level structure $\overline{\eta}$ of type $K^p$.
\end{altitemize}
These are subject to the following conditions.
\begin{altenumerate}
\item[{\rm (i)}] The Rosati involution induced by $\lambda$ is compatible with $\ast$ on $\mathcal{O}_B$.
\item[{\rm (ii)}] The determinant condition holds true.
\end{altenumerate}
In particular, if $p$ is locally nilpotent on $S$, then the associated $p$-divisible group with extra structure $(A[p^\infty],\iota|_{A[p^\infty]},\lambda|_{A[p^\infty]},\mathbb{Z}_p)$ is a $p$-divisible group with $\mathcal{D}$-structure in the sense of Definition \ref{DefPDivGroupDStr}.

Two triples $(A,\iota,\lambda,\overline{\eta})$, $(A^\prime,\iota^\prime,\lambda^\prime,\overline{\eta}^\prime)$ are said to be isomorphic if there is a prime-to-$p$-isogeny $\alpha: A\rightarrow A^\prime$ carrying $\iota$ into $\iota^\prime$, $\overline{\eta}$ into $\overline{\eta}^\prime$ and carrying $\lambda$ into a $\mathbb{Z}_{(p)}^\times$-multiple of $\lambda^\prime$ locally on $S$.
\end{definition}

We refer to \cite{KottwitzPoints}, page 390 -- 391, for the notion of a level structure of type $K^p$. Note that there is an obvious action of $\mathbf{G}(\mathbb{A}_f^p)$ by correspondences on the tower of these moduli problems.

\begin{thm} The functor $\mathfrak{M}_{K^p}$ is represented by a quasiprojective scheme $\mathcal{M}_{K^p}$ over $\mathcal{O}_{E_{\mathfrak{p}}}$. There is an isomorphism
\[
\mathcal{M}_{K^p}\otimes E_{\mathfrak{p}}\cong \bigsqcup_{\mathrm{ker}^1(\mathbb{Q},\mathbf{G})} \mathrm{Sh}_{\mathbf{G}(\mathbb{Z}_p)K^p}\otimes_E E_{\mathfrak{p}}
\]
compatible with the action of the Hecke correspondences. Here $\mathrm{ker}^1(\mathbb{Q},\mathbf{G})\subset H^1(\mathbb{Q},\mathbf{G})$ is the subset of those cohomology classes that map trivially to $H^1(\mathbb{Q}_v,\mathbf{G})$ for all places $v$ of $\mathbb{Q}$; this is a finite set.

If $C$ is a division algebra, then $\mathcal{M}_{K^p}$ is a projective variety over $\mathcal{O}_{E_{\mathfrak{p}}}$.
\end{thm}

\begin{proof} Cf. \cite{KottwitzPoints}, Section 3. We note that Kottwitz a priori only proves that these spaces are (possibly infinite) disjoint unions of quasiprojective schemes. In the generic fibre, the given description shows that there only finitely many connected components. However, the moduli spaces need not be (topologically) flat, cf. below, so that a priori there might be many junk components in the special fibre.

The problem with Kottwitz' argument is that one works with abelian varieties up to prime-to-$p$-isogeny, and the endomorphisms only live in $\mathrm{End}(A)\otimes \mathbb{Z}_{(p)}$. For any abelian scheme $A/S$, the functor $T\mapsto \mathrm{End}(A\times_S T)\otimes \mathbb{Z}_{(p)}$ on schemes over $S$ is representable by an infinite disjoint union of projective varieties over $S$. In order to restrict to finitely many components, one has to work with actual morphisms of abelian varieties whose degree is bounded.

We sketch how this can be accomplished. Choose a finitely generated $\mathbb{Z}$-lattice $M$ in $V$, look at the order $\mathcal{O}_{B,M}$ of $B$ mapping $M$ into itself, and let $K^p$ be contained in those automorphisms that fix $M\otimes \hat{\mathbb{Z}}^p\subset V\otimes \mathbb{A}_f^p$. In that case, the datum of $\overline{\eta}$ gives an actual abelian variety $A$ inside the prime-to-$p$-isogeny class, and one checks that $\mathcal{O}_{B,M}$ acts on $A$ by actual morphisms, and that a bounded multiple of $\lambda$ is an actual polarization of bounded degree. Fixing finitely many algebra generators of $\mathcal{O}_{B,M}$, one sees that each of them has bounded degree, which gives the desired statement.
\end{proof}

Moreover, for any compact open subgroup $K_p\subset \mathbf{G}(\mathbb{Z}_p)$, we introduce the cover $\mathcal{M}_{K_p,K^p}$ of $\mathcal{M}_{K^p}\otimes E_{\mathfrak{p}}$, parametrizing $K_p$-orbits of isomorphisms between $\Lambda\otimes \mathbb{Z}_p$ and the $p$-adic Tate module $T_p A$ of $A$, compatible with the $\mathcal{O}_B$-action and the hermitian forms up to a scalar. Obviously, we get an action of $\mathbf{G}(\mathbb{Z}_p)\times \mathbf{G}(\mathbb{A}_f^p)$ by correspondences on the tower of these varieties; we do not care about enlarging this action to $\mathbf{G}(\mathbb{A}_f)$ here.

\begin{prop} The cover $\pi_{K_p,K^p}: \mathcal{M}_{K_p,K^p}\rightarrow \mathcal{M}_{K^p}\otimes E_{\mathfrak{p}}$ is finite \'{e}tale, and Galois with Galois group $\mathbf{G}(\mathbb{Z}_p)/K_p$ if $K_p\subset \mathbf{G}(\mathbb{Z}_p)$ is normal. There are isomorphisms
\[
\mathcal{M}_{K_p,K^p}\cong  \bigsqcup_{\mathrm{ker}^1(\mathbb{Q},\mathbf{G})} \mathrm{Sh}_{K_pK^p}\otimes_E E_{\mathfrak{p}}
\]
compatible with the Hecke correspondences and the maps to $\mathcal{M}_{K^p}\otimes E_{\mathfrak{p}}$.
\end{prop}

\begin{proof} Once again, the existence of such isomorphisms over geometric points follows from Lemma \ref{CrucialLemma}; note that the rational $\ell$-adic Tate modules $V_\ell A$ are isomorphic to $V\otimes \mathbb{Q}_\ell$ by existence of $K^p$-level structures; hence so are $V_p A$ and $V\otimes \mathbb{Q}_p$, as the characters of $V_p A$ and $V_\ell A$ as $B$-representations agree. From here, the usual arguments imply the proposition.
\end{proof}

We have the following comparison with the local theory.

\begin{prop}\label{SerreTate} Let $\kappa$ be perfect field of characteristic $p$ over $\mathcal{O}_{E_{\mathfrak{p}}}$, and let $x\in \mathcal{M}_{K^p}(\kappa)$. Let $\underline{\overline{H}}$ be the associated $p$-divisible group with $\mathcal{D}$-structure over $\kappa$. Let $k$ be the complete unramified extension of $E_{\mathfrak{p}}$ with residue field $\kappa$. Then the complete local ring $\hat{\mathcal{O}}_{\mathcal{M}_{K^p},x}$ is isomorphic to the deformation ring $R_{\underline{\overline{H}}}$. This identifies the tubular neighborhood of $x$ in $\mathcal{M}_{K^p}^{\mathrm{rig}}\otimes_{E_{\mathfrak{p}}} k$ with $X_{\underline{\overline{H}}}$. Moreover, for any $K_p\subset \mathbf{G}(\mathbb{Z}_p)$, we have a pullback diagram
\[\xymatrix{
X_{\underline{\overline{H}},K_p}\ar[r]\ar@{^(->}[d] & X_{\underline{\overline{H}}}\ar@{^(->}[d]\\
\mathcal{M}_{K_p,K^p}^{\mathrm{rig}}\otimes_{E_{\mathfrak{p}}} k\ar[r] & \mathcal{M}_{K^p}^{\mathrm{rig}}\otimes_{E_{\mathfrak{p}}} k
}\]
In particular, $\underline{\overline{H}}$ has controlled cohomology, and for all $i\in \mathbb{Z}$, we have a $\mathrm{Gal}(\bar{k}/k)$-equivariant isomorphism
\[
(R^i\psi \pi_{K_p,K^p \ast} \mathbb{Q}_\ell)_{\overline{x}}\cong H^i(X_{\underline{\overline{H}},K_p}\otimes_k \hat{\bar{k}},\mathbb{Q}_\ell)\ .
\]
\end{prop}

\begin{proof} The first statement is a direct consequence of the Serre-Tate theorem, and the rest follows easily, using Theorem \ref{AlgebraizationControlled} for the last statements.
\end{proof}

We need local systems on the Shimura varieties. For this purpose, let $\xi$ be a finite-dimensional algebraic representation of $\mathbf{G}$ defined over a number field $L$, and let $\lambda$ be a place of $L$ above $\ell\neq p$. The usual construction, cf. \cite{KottwitzPoints}, Section 6, produces $\ell$-adic local systems $\mathcal{F}_{\xi,K^p}$, resp. $\mathcal{F}_{\xi,K_p,K^p}$ on $\mathcal{M}_{K^p}$, resp. $\mathcal{M}_{K_p,K^p}$, to which the action of the Hecke correspondences extend. We will need the following proposition.

\begin{prop}\label{CoeffDontMatter} Consider the projection $\pi_{K_p,K^p}: \mathcal{M}_{K_p,K^p}\rightarrow \mathcal{M}_{K^p}\otimes E_{\mathfrak{p}}$. Then we have a canonical isomorphism
\[
R\psi \pi_{K_p,K^p \ast} \mathcal{F}_{\xi,K_p,K^p} \cong \mathcal{F}_{\xi,K^p}\otimes R\psi \pi_{K_p,K^p \ast} \mathbb{Q}_\ell\ .
\]
\end{prop}
 
\begin{proof} Indeed,
\[\begin{aligned}
R\psi \pi_{K_p,K^p \ast} \mathcal{F}_{\xi,K_p,K^p}&\cong R\psi \pi_{K_p,K^p \ast} \pi_{K_p,K^p}^\ast \mathcal{F}_{\xi,K^p}\\
&\cong R\psi ( \mathcal{F}_{\xi,K^p}\otimes \pi_{K_p,K^p \ast} \pi_{K_p,K^p}^\ast \mathbb{Q}_\ell)\\
&\cong \mathcal{F}_{\xi,K^p} \otimes R\psi \pi_{K_p,K^p \ast} \pi_{K_p,K^p}^\ast \mathbb{Q}_\ell\\
&\cong \mathcal{F}_{\xi,K^p}\otimes R\psi \pi_{K_p,K^p \ast} \mathbb{Q}_\ell\ .
\end{aligned}\]
\end{proof}

In particular, we can define the cohomology of the Shimura variety with coefficients in the local system $\mathcal{F}_{\xi}$,
\[
H^\ast_\xi = \varinjlim_{K_p,K^p} H^\ast(\mathcal{M}_{K_p,K^p}\otimes \bar{\mathbb{Q}}_p,\mathcal{F}_{\xi,K_p,K^p})\ .
\]
This carries commuting left actions of $G_{E_{\mathfrak{p}}} = \mathrm{Gal}(\bar{\mathbb{Q}}_p/E_{\mathfrak{p}})$ and $\mathbf{G}(\mathbb{Z}_p)\times \mathbf{G}(\mathbb{A}_f^p)$. Note that in fact
\[
H^\ast_\xi = \bigoplus_{\mathrm{ker}^1(\mathbb{Q},\mathbf{G})} H^\ast_{\mathrm{Sh},\xi}
\]
with the obvious definition of $H^\ast_{\mathrm{Sh},\xi}$, compatible with all actions. The right-hand side even carries an action of $\mathbf{G}(\mathbb{A}_f)$, which we will not need, however.

Let $I_{E_{\mathfrak{p}}}\subset W_{E_{\mathfrak{p}}}\subset G_{E_{\mathfrak{p}}}$ be the inertia and Weil subgroup, and fix a geometric Frobenius element $\mathrm{Frob}\in W_{E_{\mathfrak{p}}}$. In order to formulate our main theorem, we need to introduce the notion of a Kottwitz triple.

\begin{definition} Let $j\geq 1$. Set $r:=j [\kappa_{E_{\mathfrak{p}}}:\mathbb{F}_p]$, where $\kappa_{E_{\mathfrak{p}}}$ is the residue field of $E_{\mathfrak{p}}$. A degree-$j$-Kottwitz triple $(\gamma_0;\gamma,\delta)$ consists of
\begin{altitemize}
\item a semisimple stable conjugacy class $\gamma_0\in \mathbf{G}(\mathbb{Q})$,
\item a conjugacy class $\gamma\in \mathbf{G}(\mathbb{A}_f^p)$ that is stably conjugate to $\gamma_0$, and
\item a $\sigma$-conjugacy class $\delta\in \mathbf{G}(\mathbb{Q}_{p^r})$ such that $N\delta$ is stably conjugate to $\gamma_0$
\end{altitemize}
satisfying
\begin{altenumerate}
\item[{\rm (i)}] $\gamma_0$ is elliptic in $\mathbf{G}(\mathbb{R})$,
\item[{\rm (ii)}] $\kappa_{\mathbf{G}\otimes \mathbb{Q}_p}(p\delta)=\mu^\sharp$ in $X^\ast(Z(\hat{\mathbf{G}})^{\Gamma(p)})$, where $\Gamma(p)$ is the absolute Galois group of $\mathbb{Q}_p$.
\end{altenumerate}
\end{definition}

Let $I_0$ be the centralizer of $\gamma_0$ in $\mathbf{G}$. Then Kottwitz defines a finite group $\mathfrak{K}(I_0/\mathbb{Q})$ whose Pontrjagin dual we denote by $\mathfrak{K}(I_0/\mathbb{Q})^D$. Moreover, he associates an invariant $\alpha(\gamma_0;\gamma,\delta)\in \mathfrak{K}(I_0/\mathbb{Q})^D$ to any degree-$j$-Kottwitz triple.

Finally, we can formulate our main theorem.

\begin{thm}\label{MainTheorem} Assume that the flat closure of $\mathcal{M}_{K^p}\otimes E_{\mathfrak{p}}$ in $\mathcal{M}_{K^p}$ is proper (for one and hence every $K^p$); for example, assume that $C$ is a division algebra. Let $f^p\in C_c^{\infty}(\mathbf{G}(\mathbb{A}_f^p))$, $h\in C_c^{\infty}(\mathbf{G}(\mathbb{Z}_p))$ and $\tau\in \mathrm{Frob}^j I_{E_{\mathfrak{p}}}\subset W_{E_{\mathfrak{p}}}$. Then
\[
\mathrm{tr}(\tau\times h f^p|H^{\ast}_{\mathrm{Sh},\xi}) = \sum_{\substack{(\gamma_0;\gamma,\delta)\\ \alpha(\gamma_0;\gamma,\delta)=1}} c(\gamma_0;\gamma,\delta) O_{\gamma}(f^p) TO_{\delta\sigma}(\phi_{\tau,h}) \mathrm{tr}\ \xi(\gamma_0)\ ,
\]
where the sum runs over degree-$j$-Kottwitz triples, and $c(\gamma_0;\gamma,\delta)$ is a volume factor defined as in \cite{KottwitzPoints}, p. 441. The Haar measures on $\mathbf{G}(\mathbb{Q}_p)$ resp. $\mathbf{G}(\mathbb{Q}_{p^r})$ are normalized by giving $\mathbf{G}(\mathbb{Z}_p)$ resp. $\mathbf{G}(\mathbb{Z}_{p^r})$ volume $1$.
\end{thm}

Slightly more generally, we will prove this theorem for any function $\phi_{\tau,h}^\prime\in C_c^\infty(\mathbf{G}(\mathbb{Q}_{p^r}))$ in place of $\phi_{\tau,h}$ that has the following properties:
\begin{altenumerate}
\item[{\rm (i)}] $\phi_{\tau,h}^\prime(\delta)=0$ unless $\delta$ is associated to some $p$-divisible group with $\mathcal{D}$-structure $\underline{\overline{H}}$;
\item[{\rm (ii)}] If $\underline{\overline{H}}_{\bar{\mathbb{F}}_p}$ is the $p$-divisible group with $\mathcal{D}$-structure associated to some point of $\mathcal{M}_{K^p}(\bar{\mathbb{F}}_p)$ (for some $K^p$), then $\phi_{\tau,h}^\prime(\delta)=\phi_{\tau,h}(\delta)$.
\end{altenumerate}

Conjecturally, one would expect that any $p$-divisible group with $\mathcal{D}$-structure occurs in the Shimura variety as in (ii), and hence that $\phi_{\tau,h}^\prime$ is unique. In \cite{ScholzeShin}, we will use this argument the other way around: We will show that (the twisted orbital integrals of) $\phi_{\tau,h}$ are uniquely determined by the formula of the theorem in some cases, and deduce that every $p$-divisible group with $\mathcal{D}$-structure occurs in the Shimura variety.

We note that it may happen that the models $\mathcal{M}_{K^p}$ are not flat, not even topologically flat, thereby explaining the formulation of the properness assumption in the theorem above. In general, such questions are the subject of the theory of local models of Shimura varieties. We only mention here that the results of Pappas and Rapoport in \cite{PappasRapoportLocalModels1} show that in case A, the model that we have given, also known as the naive local model, need not be topologically flat. In fact, the dimension of the special fibre may be larger than the dimension of the generic fibre. However, in case C, results of G\"ortz, \cite{GoertzLocalModels}, show that the model is always topologically flat.

In \cite{ScholzeShin}, we will need another case where we can apply Theorem \ref{MainTheorem}. This relies on a theorem of K.-W. Lan, \cite{Lan}. Recall that in case A, the group $\mathbf{G}$ is a unitary similitude group sitting in an exact sequence
\[
0\rightarrow \mathrm{Res}_{F_0/\mathbb{Q}} \mathbf{G}_1\rightarrow \mathbf{G}\rightarrow \mathbb{G}_m\rightarrow 0\ .
\]

\begin{thm} Assume that in case A, the group $\mathbf{G}_1/F_0$ is compact at one infinite place of $F_0$, i.e. is isomorphic to $U(0,n)$. Then the flat closure of $\mathcal{M}_{K^p}\otimes E_{\mathfrak{p}}$ in $\mathcal{M}_{K^p}$ is proper.
\end{thm}

\begin{rem} The proof does not show that $\mathcal{M}_{K^p}$ itself is proper.
\end{rem}

\begin{proof} We use Theorem 5.3.3.1 of \cite{Lan}. We feel that it is worthwhile to make explicit the objects occuring in this theorem. First, the set of primes, denoted $\Box$ in \cite{Lan}, can be chosen to be the empty set. Then $M_{\mathcal{H}}$ lives over $S_0=\mathrm{Spec}\, E$ and becomes after base-change to $E_{\mathfrak{p}}$ equal to $\mathcal{M}_{K^p}\otimes E_{\mathfrak{p}}$, where $\mathcal{H}$ and $K^p$ are supposed to give corresponding level structures. We take $M^\prime$ equal to the flat closure of $\mathcal{M}_{K^p}\otimes E_{\mathfrak{p}}$ in $\mathcal{M}_{K^p}$, living over $S^\prime = \mathrm{Spec}\, \mathcal{O}_{E_{\mathfrak{p}}}$. Let $S_1^\prime = \mathrm{Spec}\, E_{\mathfrak{p}}$. Then condition (1) is satisfied tautologically, and condition (2) is verified as in \cite{KottwitzPoints}, page 392. Finally, condition (3) is ensured by the signature condition, cf. Remark 5.3.3.2 in \cite{Lan}.
\end{proof}

\section{Fixed points of correspondences}

In order to prove the theorem, we can make the following assumptions. First, fix a sufficiently small compact open subgroup $K^p\subset \mathbf{G}(\mathbb{A}_f^p)$ such that $f^p$ is bi-$K^p$-invariant. In fact, assume that $f^p$ is the characteristic function of $K^pg^pK^p$ divided by the volume of $K^p$ for some $g^p\in \mathbf{G}(\mathbb{A}_f^p)$. Further, let $K_p\subset \mathbf{G}(\mathbb{Z}_p)$ be a normal subgroup, such that $h$ is bi-$K_p$-invariant; we assume that $h$ is the characteristic function of $K_pg_p$ divided by the volume of $K_p$ for some $g_p\in \mathbf{G}(\mathbb{Z}_p)$. We have the following diagram, where $K^p_{g^p} = K^p\cap (g^p)^{-1}K^pg^p$.
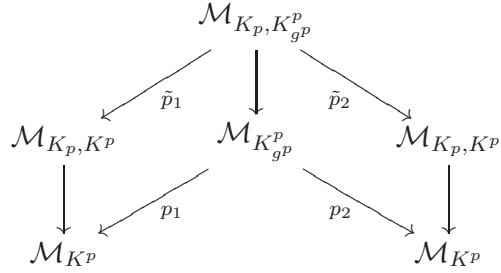
\[\xymatrix{
& \mathcal{M}_{K_p,K^p_{g^p}}\ar[ld]^{\tilde{p}_1} \ar[d] \ar[rd]_{\tilde{p}_2} & \\
\mathcal{M}_{K_p,K^p}\ar[d] & \mathcal{M}_{K^p_{g^p}}\ar[ld]^{p_1} \ar[rd]_{p_2} & \mathcal{M}_{K_p,K^p}\ar[d] \\
\mathcal{M}_{K^p} & & \mathcal{M}_{K^p}
}\]
Here, the left-hand diagonal projections are the natural projections, whereas the right-hand diagonal projections are (in the upper case) the composite of the natural projection $\mathcal{M}_{K_p,K^p_{g^p}}\rightarrow \mathcal{M}_{K_p,(g^p)^{-1}K^pg^p}$ with the action of $(g_p,g^p)$, which gives an isomorphism $\mathcal{M}_{K_p,(g^p)^{-1}K^pg^p}\cong \mathcal{M}_{K_p,K^p}$.

Let $\Phi_{\mathfrak{p}}: \mathcal{M}_{K^p}\otimes \kappa_{E_{\mathfrak{p}}}\rightarrow \mathcal{M}_{K^p}\otimes \kappa_{E_{\mathfrak{p}}}$ be the Frobenius morphism. Our first goal is to describe the fixed points of the lower correspondence restricted to the special fibre composed with the Frobenius correspondence $\Phi_{\mathfrak{p}}^j$. \footnote{Here, we follow the conventions on correspondences and associated maps on cohomology that is used in \cite{Varshavsky}. In particular, we compose $p_1$ with the Frobenius correspondence; in \cite{KottwitzPoints}, $p_2$ is composed with the Frobenius correspondence instead. We note that Kottwitz uses the same correspondence associated to $g^p$, but constructs out of it a map on cohomology that goes the other way; this results in giving the inverse of the action of $g^p$ on the cohomology, i.e. the action of $(g^p)^{-1}$. This leads to the occurence of the characteristic function of $K^p (g^p)^{-1} K^p$ later in Kottwitz' paper, instead of $f^p$.} In fact, remembering that our model is not always topologically flat, it will be enough to describe those fixed points that lie in the closure of the generic fibre. This is done (under slightly stronger, but unnecessary assumptions) in \cite{KottwitzPoints}. Let us recall the description.

By definition, a fixed point of the correspondence is a point $(\overline{A},\iota,\lambda,\overline{\eta})\in \mathcal{M}_{K^p_{g^p}}(\bar{\mathbb{F}}_p)$ such that $(\overline{A},\iota,\lambda,\overline{\eta})$ and $\sigma^r(\overline{A},\iota,\lambda,\overline{\eta g^p})$ define the same point of $\mathcal{M}_{K^p}(\bar{\mathbb{F}}_p)$, where $\sigma$ is the $p$-th power map on $\bar{\mathbb{F}}_p$ and $\sigma^r(\overline{A},\ldots)$ is obtained through extension of scalars along $\sigma^r$ from $(\overline{A},\ldots)$.

This translates into the condition that there is some prime-to-$p$-isogeny $u: \sigma^r(\overline{A})\rightarrow \overline{A}$ compatible with $\iota$ and sending $\sigma^r(\overline{\eta g^p})$ into $\overline{\eta}$, and such that there is some $c_0\in \mathbb{Z}_{(p)}^\times$ with $u^\ast \lambda = c_0 \sigma^r(\lambda)$. At this point, let us recall some definitions from \cite{KottwitzPoints}.

\begin{definition} \begin{altenumerate}
\item[{\rm (i)}] A virtual abelian variety over $\mathbb{F}_{p^r}$ is a pair $A=(\overline{A},u)$ consisting of an abelian variety $\overline{A}$ up to prime-to-$p$-isogeny over $\bar{\mathbb{F}}_p$ with a prime-to-$p$-isogeny $u: \sigma^r(\overline{A})\rightarrow \overline{A}$.
\item[{\rm (ii)}] A homomorphism between two virtual abelian varieties $A_1$, $A_2$ is a homomorphism $f: \overline{A}_1\rightarrow \overline{A}_2$ such that $f u_1 = u_2 \sigma^r(f)$.
\item[{\rm (iii)}] The Frobenius morphism $\pi_A\in \mathrm{End}(A)$ of a virtual abelian variety is the composition $u\circ \Phi^r$, where $\Phi: \overline{A}\rightarrow \sigma(\overline{A})$ is the relative Frobenius over $\bar{\mathbb{F}}_p$.
\item[{\rm (iv)}] For a rational number $c$ of the form $c=p^rc_0$ with $c_0\in \mathbb{Z}_{(p)}^\times$, a $c$-polarization of a virtual abelian $A=(\overline{A},u)$ is a $\mathbb{Q}$-polarization $\lambda: \overline{A}\rightarrow \overline{A}^{\vee}$ (i.e., $\lambda$ is only a quasi-isogeny, and some multiple of $\lambda$ is a polarization) such that $u^\ast \lambda = c_0 \sigma^r(\lambda)$.
\end{altenumerate}
\end{definition}

We remark that the category of abelian varieties over $\mathbb{F}_{p^r}$ forms a full subcategory of the category of virtual abelian varieties over $\mathbb{F}_{p^r}$ in the obvious way. Moreover, a polarization of an abelian variety over $\mathbb{F}_{p^r}$ induces a $p^r$-polarization of the associated virtual abelian variety.

\begin{prop}\label{VirtAbVarPDiv} Let $A=(\overline{A},u)$ be a virtual abelian variety over $\mathbb{F}_{p^r}$. Then $\overline{A}[p^\infty]$ descends via $u|_{\sigma^r(\overline{A}[p^\infty])}$ to a $p$-divisible group $H$ over $\mathbb{F}_{p^r}$, which we sometimes also denote by $A[p^\infty]$. Moreover, if $\lambda$ is a $c$-polarization of $A$, and $d\in \mathbb{Z}_{p^r}$ has norm $Nd = dd^\sigma\cdots d^{\sigma^{r-1}}=c_0^{-1}$, and $\mathbb{L}$ is the $1$-dimensional $\mathbb{Z}_p$-local system over $\mathbb{F}_{p^r}$ corresponding to $d$, then $\lambda$ induces a twisted principal polarization
\[
\lambda[p^\infty]: A[p^\infty]\rightarrow A[p^\infty]^\vee\otimes \mathbb{L}\ .
\]
\end{prop}

\begin{proof} As $u$ is a prime-to-$p$-isogeny, it becomes an isomorphism on the $p$-divisible group. To check that the descent is effective, it suffices to check this on $p^m$-torsion points for all $m$; but restricted to the $p^m$-torsion points, everything is defined over a finite subextension of $\bar{\mathbb{F}}_p$. Under the corresponding identification of $\sigma^k(\overline{A}[p^m])$ with $\overline{A}[p^m]$ for $k$ large, some power of $u$ has to be the identity morphism (because there are only finitely many automorphisms of a finite flat group scheme over a finite field, as the ring of global sections is finite). This shows that we actually get a Galois descent datum. We leave the verification about polarizations to the reader.
\end{proof}

In particular, from a fixed point of the correspondence, we get a $c$-polarized virtual abelian variety $A$ over $\mathbb{F}_{p^r}$, with an action $\iota$ of $\mathcal{O}_B$, compatible with the polarization $\lambda$, i.e. a triple $(A,\iota,\lambda)$. Note that the previous proposition associates to $(A,\iota,\lambda)$ a $p$-divisible group with $\mathcal{D}$-structure $\underline{\overline{H}} = (A[p^\infty],\iota|_{A[p^\infty]},\lambda[p^\infty],\mathbb{L})$ over $\mathbb{F}_{p^r}$. Conversely a triple $(A,\iota,\lambda)$ whose associated $p$-divisible group with extra structure is a $p$-divisible group with $\mathcal{D}$-structure, together with a suitable level structure $\overline{\eta}$ of type $K^p_{g^p}$, will give a fixed point of the correspondence. 

By the existence of a level structure of type $K^p$, we know that for all $\ell\neq p$, the rational $\ell$-adic Tate module $V_\ell \overline{A}$ is isomorphic to $V\otimes \mathbb{Q}_\ell$; fixing an isomorphism, the Frobenius morphism $\pi_A\in \mathrm{End}(A)$ gives rise to a $B$-linear automorphism of $V_\ell$; we define $\gamma_\ell\in \mathbf{G}(\mathbb{Q}_\ell)$ as its inverse. Its conjugacy class is well-defined, and the elements combine into a conjugacy class $\gamma\in \mathbf{G}(\mathbb{A}_f^p)$.

Now assume that $x\in \mathcal{M}_{K^p_{g^p}}(\bar{\mathbb{F}}_p)$ lies in the closure of the generic fibre, so that $X_{\underline{\overline{H}}}\neq \emptyset$. Then our local considerations give an element $\delta\in \mathbf{G}(\mathbb{Q}_{p^r})$, well-defined up to $\sigma$-conjugation by $\mathbf{G}(\mathbb{Z}_{p^r})$. It satisfies $\kappa_{\mathbf{G}\otimes \mathbb{Q}_p} (p\delta) = \mu^\sharp$. We have the following proposition.

\begin{prop} There is a unique semisimple stable conjugacy class $\gamma_0\in \mathbf{G}(\mathbb{Q})$ such that $(\gamma_0;\gamma,\delta)$ is a degree-$j$-Kottwitz triple.
\end{prop}

\begin{proof} Uniqueness is clear; for existence, follow the arguments in \cite{KottwitzPoints}, Section 14.
\end{proof}

\begin{prop} The invariant $\alpha(\gamma_0;\gamma,\delta)\in \mathfrak{K}(I_0/\mathbb{Q})^D$ is trivial.
\end{prop}

\begin{proof} Follow the arguments in \cite{KottwitzPoints}, Section 15.
\end{proof}

Our next aim is to parametrize fixed points in the isogeny class of the virtual abelian variety with extra structure $(A,\iota,\lambda)$, following the discussion in \cite{KottwitzPoints}, page 431 -- 433. Let $I/\mathbb{Q}$ be the group of self-quasiisogenies of $(A,\iota,\lambda)$ (preserving $\lambda$ up to a scalar in $\mathbb{Q}^\times$); it is an inner form of the centralizer $I_0$ of $\gamma_0$ in $\mathbf{G}$. Consider the set $Y$ of fixed points $y\in \mathcal{M}_{K^p}(\bar{\mathbb{F}}_p)$ that lie in the closure of the generic fibre, equipped with a $B$-linear quasiisogeny $\varphi: A^\prime\rightarrow A$ preserving the polarization up to a scalar in $\mathbb{Q}^\times$, from the associated triple $(A^\prime,\iota^\prime,\lambda^\prime)$ to $(A,\iota,\lambda)$. The set of such fixed points itself is then given by $I(\mathbb{Q})\backslash Y$. Let $N$ be the rational Dieudonn\'{e} module of $A[p^\infty]$, which is a $\mathbb{Q}_{p^r}$-vector space.

\begin{prop} There is an injection from $Y$ into the set of $\mathcal{O}_B$-stable selfdual $\mathbb{Z}_{p^r}$-lattices $\Lambda^\prime\subset N$ satisfying $p\Lambda^\prime\subset p\delta\sigma \Lambda^\prime\subset \Lambda^\prime$, together with an element $z\in \mathbf{G}(\mathbb{A}_f^p)/K^p$ satisfying $z^{-1}\gamma z\in g^p K^p$. In that case, $(\Lambda^\prime,F=p\delta\sigma)$ is the Dieudonn\'{e} module associated to $A^\prime[p^\infty]$ (with extra structure).
\end{prop}

\begin{proof} As we consider abelian varieties up to prime-to-$p$-isogeny, it is clear that $\Lambda^\prime$ determines $(A^\prime,\iota^\prime,\lambda^\prime)$. The element $z$ describes the level structure of type $K^p_{g^p}$, cf. \cite{KottwitzPoints}, page 432. We note that as we are using slightly different conventions on correspondences, a fixed point satisfies $\overline{\eta} = \sigma^r(\overline{\eta g^p})$; this translates into $\eta = \pi_A \eta g^p$ modulo $K^p$, or equivalently $z = \gamma^{-1} z g^p$ modulo $K^p$, which means that $z^{-1} \gamma z\in g^p K^p$.
\end{proof}

By Lemma \ref{CrucialLemma} (and the remark following it), the lattice $\Lambda^\prime$ can be equivalently described by giving an element $w\in \mathbf{G}(\mathbb{Q}_{p^r})/\mathbf{G}(\mathbb{Z}_{p^r})$, where $\Lambda^\prime = w\Lambda$. The corresponding $\delta^\prime\in \mathbf{G}(\mathbb{Q}_{p^r})$ up to $\sigma$-conjugation by $\mathbf{G}(\mathbb{Z}_{p^r})$ is then given by $\delta^\prime = w^{-1} \delta w^{\sigma}$.

We will check in the next section that a fixed point giving rise to $(\Lambda^\prime,z)$ contributes the summand
\[
\phi_{\tau,h}(\delta^\prime) \mathrm{tr}\ \xi(\gamma_0)
\]
to the Lefschetz trace formula. Conversely, note that if $(\Lambda^\prime,z)$ is as in the previous proposition and $\phi_{\tau,h}(\delta^\prime)\neq 0$, then necessarily the pair $(\Lambda^\prime,z)$ comes from an element of $Y$: Namely, the function $\phi_{\tau,h}$ is nonzero only at those elements whose associated $p$-divisible group with extra structure does satisfy the determinant condition, i.e. is a $p$-divisible group with $\mathcal{D}$-structure.

Summarizing, we get the following proposition.

\begin{prop} The contribution of fixed points isogenous to $(A,\iota,\lambda)$ is given by
\[
\mathrm{vol}(I(\mathbb{Q})\backslash I(\mathbb{A}_f)) O_{\gamma}(f^p) TO_{\delta\sigma}(\phi_{\tau,h}) \mathrm{tr}\ \xi(\gamma_0)\ .
\]
\end{prop}

\begin{proof} For the remaining easy verifications, see \cite{KottwitzPoints}, page 432.
\end{proof}

Our main theorem now follows from the results of the application of the Lefschetz trace formula given in the next section (affirming the form of the contribution of fixed points used above) by going through the rest of the paper of Kottwitz, Sections 17 -- 19, which give a full description of the set of isogeny classes.     

\section{Application of the Lefschetz trace formula}

We want to evaluate $\mathrm{tr}(\tau\times hf^p|H^{\ast}_{\xi})$ via the Lefschetz trace formula. Let us recall the correspondence
\[\xymatrix{
& \mathcal{M}_{K_p,K^p_{g^p}}\ar[ld]^{\tilde{p}_1} \ar[d] \ar[rd]_{\tilde{p}_2} & \\
\mathcal{M}_{K_p,K^p}\ar[d] & \mathcal{M}_{K^p_{g^p}}\ar[ld]^{p_1} \ar[rd]_{p_2} & \mathcal{M}_{K_p,K^p}\ar[d] \\
\mathcal{M}_{K^p} & & \mathcal{M}_{K^p}
}\]
The upper correspondence in the diagram above extends canonically to a cohomological correspondence $u: \tilde{p}_{2!}\tilde{p}_1^{\ast}\mathcal{F}_{\xi,K_p,K^p}\longrightarrow \mathcal{F}_{\xi,K_p,K^p}$ induced from the action of $(g_p,g^p)$. Let
\[
(g_p,g^p)_\ast:  H^\ast(\mathcal{M}_{K_p,K^p}\otimes \bar{\mathbb{Q}}_p,\mathcal{F}_{\xi,K_p,K^p})\longrightarrow H^\ast(\mathcal{M}_{K_p,K^p}\otimes \bar{\mathbb{Q}}_p,\mathcal{F}_{\xi,K_p,K^p})
\]
be the associated map on cohomology. Of course, $\tau$ also acts on the cohomology, and it is a standard fact that
\[
\mathrm{tr}(\tau\times hf^p|H^{\ast}_{\xi}) = \tr(\tau\times (g_p,g^p)_\ast|H^\ast(\mathcal{M}_{K_p,K^p}\otimes \bar{\mathbb{Q}}_p,\mathcal{F}_{\xi,K_p,K^p}))\ .
\]

We rewrite the cohomology group using proper base change and Proposition \ref{CoeffDontMatter}
\[\begin{aligned}
H^\ast(\mathcal{M}_{K_pK^p}\otimes \bar{\mathbb{Q}}_p,\mathcal{F}_{\xi,K_p,K^p}) &= H^\ast(\mathcal{M}_{K^p}\otimes \bar{\mathbb{Q}}_p,\pi_{K_p,K^p \ast}\mathcal{F}_{\xi,K_p,K^p})\\
&= H^\ast(\mathcal{M}_{K^p}\otimes \bar{\mathbb{F}}_p,R\psi\pi_{K_p,K^p \ast}\mathcal{F}_{\xi,K_p,K^p})\\
&= H^\ast(\mathcal{M}_{K^p}\otimes \bar{\mathbb{F}}_p, \mathcal{F}_{\xi,K^p}\otimes R\psi \pi_{K_p,K^p \ast} \mathbb{Q}_\ell)\ .
\end{aligned}\]
We have cohomological correspondences $p_{2!} p_1^\ast \mathcal{F}_{\xi,K^p}\rightarrow \mathcal{F}_{\xi,K^p}$ induced from $g^p$ and $p_{2!} p_1^\ast R\psi \pi_{K_p,K^p \ast} \mathbb{Q}_\ell\rightarrow R\psi \pi_{K_p,K^p \ast} \mathbb{Q}_\ell$ induced from $g_p$. Their tensor product gives a cohomological correspondence on $\mathcal{F}_{\xi,K^p}\otimes R\psi \pi_{K_p,K^p \ast} \mathbb{Q}_\ell$ which induces in cohomology the map $(g_p,g^p)_\ast$ under the isomorphism of cohomology groups
\[
H^\ast(\mathcal{M}_{K_pK^p}\otimes \bar{\mathbb{Q}}_p,\mathcal{F}_{\xi,K_p,K^p}) \cong H^\ast(\mathcal{M}_{K^p}\otimes \bar{\mathbb{F}}_p, \mathcal{F}_{\xi,K^p}\otimes R\psi \pi_{K_p,K^p \ast} \mathbb{Q}_\ell)\ .
\]

Finally, because $\mathcal{M}_{K^p}\otimes \bar{\mathbb{F}}_p$ is defined over the finite field $\kappa_{E_{\mathfrak{p}}}$, we also have the Frobenius correspondence
\[
\mathcal{M}_{K^p}\otimes \bar{\mathbb{F}}_p\buildrel {\Phi_{\mathfrak{p}}^j}\over \longleftarrow \mathcal{M}_{K^p}\otimes \bar{\mathbb{F}}_p\buildrel {=}\over\rightarrow \mathcal{M}_{K^p}\otimes \bar{\mathbb{F}}_p\ .
\]
The action of $W_{E_{\mathfrak{p}}}$ on the vanishing cycles gives rise to a cohomological correspondence $\Phi_{\mathfrak{p}}^{j \ast} R\psi \pi_{K_p,K^p \ast} \mathbb{Q}_\ell\rightarrow R\psi \pi_{K_p,K^p \ast} \mathbb{Q}_\ell$ induced by $\tau\in \mathrm{Frob}^j I_{E_{\mathfrak{p}}}$, and because $\mathcal{F}_{\xi,K^p}$ is also defined over $\kappa_{E_{\mathfrak{p}}}$, we have a cohomological correspondence $\Phi_{\mathfrak{p}}^{j \ast} \mathcal{F}_{\xi,K^p}\rightarrow \mathcal{F}_{\xi,K^p}$. Their tensor product gives a cohomological correspondence on $\mathcal{F}_{\xi,K^p}\otimes R\psi \pi_{K_p,K^p \ast} \mathbb{Q}_\ell$ which in cohomology
\[
H^\ast(\mathcal{M}_{K_pK^p}\otimes \bar{\mathbb{Q}}_p,\mathcal{F}_{\xi,K_p,K^p}) \cong H^\ast(\mathcal{M}_{K^p}\otimes \bar{\mathbb{F}}_p, \mathcal{F}_{\xi,K^p}\otimes R\psi \pi_{K_p,K^p \ast} \mathbb{Q}_\ell)
\]
gives the action of $\tau$.

Now we take the composite correspondence
\[
\mathcal{M}_{K^p}\otimes \bar{\mathbb{F}}_p\buildrel {\Phi_{\mathfrak{p}}^j\circ p_1}\over \longleftarrow \mathcal{M}_{K^p_{g^p}}\otimes \bar{\mathbb{F}}_p\buildrel {p_2}\over\rightarrow \mathcal{M}_{K^p}\otimes \bar{\mathbb{F}}_p
\]
with the composite cohomological correspondence
\[
u: p_{2!} (\Phi_{\mathfrak{p}}^j\circ p_1)^\ast \mathcal{F}_{\xi,K^p}\otimes R\psi \pi_{K_p,K^p \ast} \mathbb{Q}_\ell\rightarrow \mathcal{F}_{\xi,K^p}\otimes R\psi \pi_{K_p,K^p \ast} \mathbb{Q}_\ell\ .
\]
This induces a map on
\[
H^\ast(\mathcal{M}_{K_pK^p}\otimes \bar{\mathbb{Q}}_p,\mathcal{F}_{\xi,K_p,K^p}) \cong H^\ast(\mathcal{M}_{K^p}\otimes \bar{\mathbb{F}}_p, \mathcal{F}_{\xi,K^p}\otimes R\psi \pi_{K_p,K^p \ast} \mathbb{Q}_\ell)
\]
which is given by $\tau\times (g_p,g^p)_\ast$ by our previous considerations. We are interested in calculating its trace.

Now we use the Lefschetz trace formula in the form given in Theorem 2.3.2 b) of \cite{Varshavsky}.

\begin{thm} The Lefschetz trace formula gives
\[
\mathrm{tr}(\tau \times hf^p | H^\ast_\xi) = \sum_{\substack{x\in \mathcal{M}_{K^p_{g^p}}(\bar{\mathbb{F}}_p)\\ (\Phi_{\mathfrak{p}}^j\circ p_1)(x) = p_2(x)}} \mathrm{tr}(u_x)\ ,
\]
where the local term $\mathrm{tr}(u_x)$ is the naive local term given as the trace of the morphism
\[\begin{aligned}
u_x: (\mathcal{F}_{\xi,K^p}\otimes R\psi \pi_{K_p,K^p \ast} \mathbb{Q}_\ell)_{(\Phi_{\mathfrak{p}}^j\circ p_1)(x)}&\cong ((\Phi_{\mathfrak{p}}^j\circ p_1)^\ast (\mathcal{F}_{\xi,K^p}\otimes R\psi \pi_{K_p,K^p \ast} \mathbb{Q}_\ell))_x\\
&\rightarrow (p_2^\ast (\mathcal{F}_{\xi,K^p}\otimes R\psi \pi_{K_p,K^p \ast} \mathbb{Q}_\ell))_x\\
&\cong (\mathcal{F}_{\xi,K^p}\otimes R\psi \pi_{K_p,K^p \ast} \mathbb{Q}_\ell)_{p_2(x)}\ ,
\end{aligned}\]
the middle map being induced by $u$ (noting that $p_{2!}$ is left-adjoint to $p_2^\ast$, $p_2$ being \'{e}tale), and the outer two terms are identified, $x$ being a fixed point.$\hfill \Box$
\end{thm}

Note that the local term vanishes unless $x$ lies in the closure of the generic fibre; in that case we have associated a degree-$j$-Kottwitz triple $(\gamma_0;\gamma,\delta)$ to $x$. The morphism $u_x$ naturally factors as the tensor product of its actions on $\mathcal{F}_{\xi,K^p}$ and $R\psi \pi_{K_p,K^p \ast} \mathbb{Q}_\ell$, so that $\mathrm{tr}(u_x)$ factors into the product of the corresponding terms. The trace on $\mathcal{F}_{\xi,K^p}$ is computed by Kottwitz to be $\mathrm{tr}\ \xi(\gamma_0)$ in \cite{KottwitzPoints}, pages 433 -- 434. We are left with the trace on $R\psi \pi_{K_p,K^p \ast} \mathbb{Q}_\ell$. Using Proposition \ref{SerreTate}, one may express the vanishing cycles in terms of the cohomology of the deformation spaces of the $p$-divisible group $\underline{\overline{H}}$ over $\mathbb{F}_{p^r}$ associated to the fixed point $x$ as in the last section. Going through all identifications, this shows that the second factor is exactly $\phi_{\tau,h}(\delta)$, as desired.

\bibliographystyle{abbrv}
\bibliography{DeformationSpaces}

\end{document}